\documentclass{article}

\usepackage{PRIMEarxiv}

\usepackage[utf8]{inputenc} 
\usepackage[T1]{fontenc}    
\usepackage{hyperref}       
\usepackage{url}            
\usepackage{booktabs}       
\usepackage{amsfonts}       
\usepackage{amsmath}
\usepackage{nicefrac}       
\usepackage{microtype}      
\usepackage{lipsum}
\usepackage{graphicx}
\usepackage{subfigure}
\usepackage{mathpazo}
\usepackage{sectsty}
\usepackage{tikz}

\graphicspath{{media/}}     

\def\bfa{{\mathbf{a}}}
\def\bfb{{\mathbf{b}}}

\def\bfg{{\mathbf{g}}}

\def\bfu{{\mathbf{u}}}
\def\bfv{{\mathbf{v}}}
\def\bfx{{\mathbf{x}}}

\def\bfz{{\mathbf{z}}}
\def\bfdelta{{\boldsymbol{\delta}}}

\def\bfphi{{\boldsymbol{\phi}}}

\def\bfxi{{\boldsymbol{\xi}}}

\def\du{\textrm{d}\bfu}

\def\dxi{\textrm{d}\bfxi}

\def\calG{\mathcal{G}}

\def\calU{\mathcal{U}}

\def\calJ{\mathcal{J}}
\def\calL{\mathcal{L}}

\newcommand{\change}[1]{\textcolor{black}{#1}}
  
\title{A Survey on Design-space Dimensionality Reduction Methods for Shape Optimization
}

\author{
  Andrea Serani$^{1,\star}$ and Matteo Diez$^1$\\
  $^1$National Research Council-Institute of Marine Engineering, Rome, Italy\\
  $^\star$\texttt{andrea.serani@cnr.it} \\
}

\begin{document}

\begin{tikzpicture}[remember picture,overlay]
   \node [rectangle, fill=cyan, fill opacity=0.5, anchor=north, minimum width=\paperwidth, minimum height=3cm] at (current page.north) {};

   \node [anchor=north, minimum width=\paperwidth, minimum height=3cm, text width=\textwidth, align=center, text height=5ex, text depth=15ex, align=left] at (current page.north) {
     \sffamily\small
     \textbf{This is a preprint submitted to:} \textit{Archives of Computational Methods in Engineering}
   };
\end{tikzpicture}

\maketitle

\begin{abstract}
The rapidly evolving field of engineering design of functional surfaces necessitates sophisticated tools to manage the inherent complexity of high-dimensional design spaces. This \change{survey paper offers a scoping review, i.e., a literature mapping synthesis borrowed from clinical medicine, delving} into the field of design-space dimensionality reduction techniques tailored for shape optimization, bridging traditional methods and cutting-edge technologies. Dissecting the spectrum of these techniques, from classical linear approaches like principal component analysis to more nuanced nonlinear methods such as autoencoders, the discussion extends to innovative physics-informed methods that integrate physical data into the dimensionality reduction process, enhancing \change{the physical relevance and effectiveness of reduced design spaces}. By integrating these methods into optimization frameworks, it is shown how they significantly mitigate the curse of dimensionality, streamline computational processes, and refine the \change{design} exploration and optimization of complex functional surfaces. The survey provides a classification of method and highlights the transformative impact of these techniques in simplifying design challenges, thereby fostering more efficient and effective engineering solutions.
\end{abstract}

\keywords{Shape optimization \and Space reduction \and Dimensionality reduction \and Representation learning \and  Principal component analysis \and Parametric model embedding \and Active subspace \and Autoencoders \and Scoping review}

\section{Introduction}
\begin{figure*}[!b]
\centering
\includegraphics[width=0.75\textwidth]{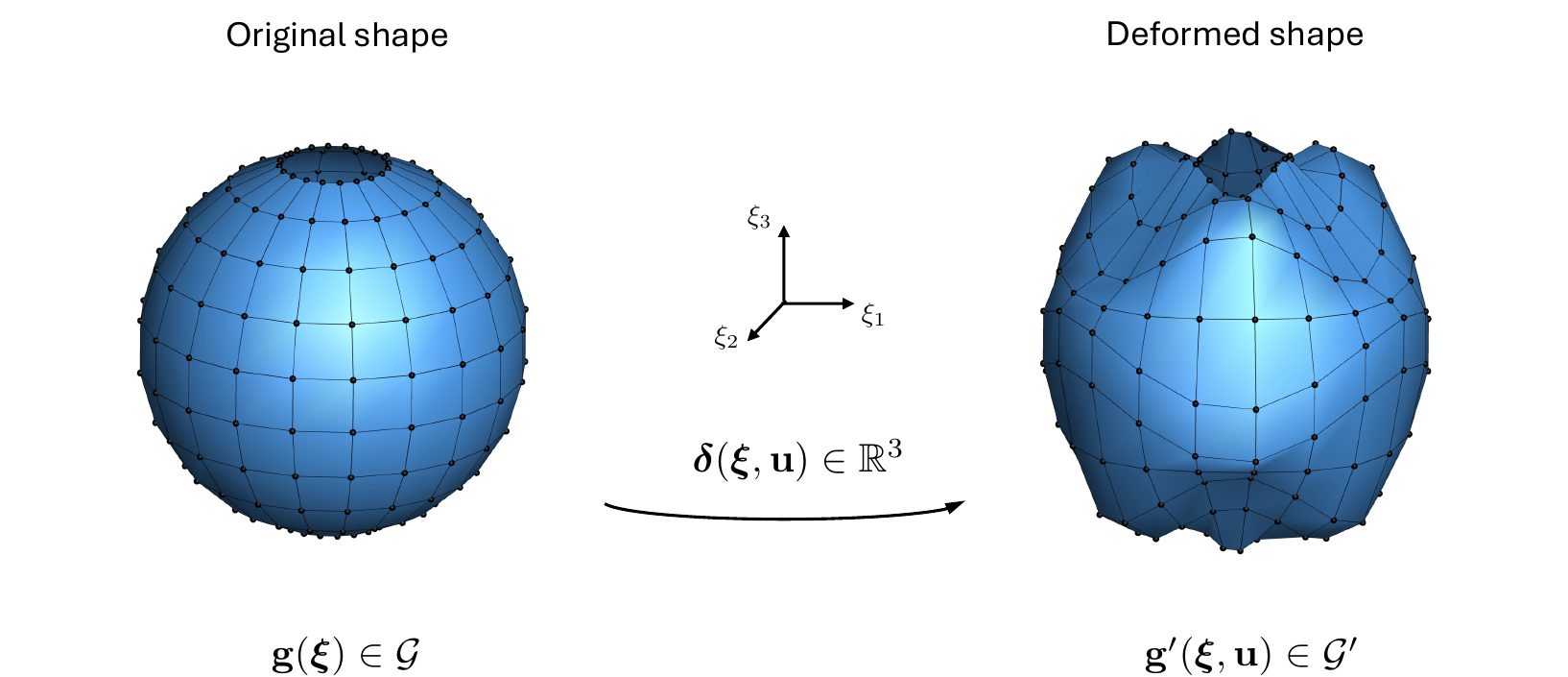}
\caption{Shape modification example}\label{fig:deform}
\end{figure*}

Shape optimization is integral to advancing performance and efficiency in a myriad of engineering disciplines, including aerodynamics, hydrodynamics, and heat transfer, \change{with applications ranging from turbomachinery to aircraft design, from ships and underwater vehicles to cars and motorcycle, from energy harvesting devices to space rockets}. The optimization process is crucial not only for enhancing functional outputs but also for meeting increasingly stringent environmental regulations and safety standards. Each of these fields involves navigating a complex, high-dimensional design space defined by an extensive array of design variables and constraints, which complicates the optimization process and amplifies computational demands.

The challenge is further intensified by the so-called \emph{curse of dimensionality} \cite{bellman1966dynamic}, which refers to the exponential growth in computational resources required as the number of dimensions (variables) in the problem increases. 
To manage this complexity, dimensionality reduction techniques have emerged as a cornerstone strategy. By effectively reducing the number of variables under consideration without significant loss of critical information, these techniques simplify the design space, making it more tractable for optimization.

\change{Although dimensionality reduction methods are widely used in various engineering applications, shape optimization problems pose unique challenges and requirements. Unlike typical functional-space applications (e.g., flow fields, stress/strain fields), shape optimization inherently involves geometric design spaces, requiring explicit and careful adaptations of these techniques to address geometric constraints, ensure feasible physical back-projections, and embed meaningful geometric metrics. Thus, a critical analysis specifically tailored to the geometric peculiarities and demands of shape optimization is essential, as traditional dimension reduction methods cannot be applied directly without significant assumptions and methodological modifications.}

This survey paper provides a comprehensive \change{scoping} review on dimensionality reduction techniques applied in the context of shape optimization. It provides a first classification of methods and explores a variety of approaches from traditional linear to more advanced nonlinear methods. Each technique's theoretical foundations, applications, and impacts on the \change{shape} optimization process are critically examined, highlighting their integration within simulation-based design optimization frameworks. 

\change{It is worth clarifying here that, while the term \emph{shape optimization} broadly encompasses a variety of problems, including structural, topological, and multi-patch manifold optimization, this survey specifically addresses the optimization of external shapes. By definition, these external shapes represent the primary geometric interfaces interacting directly with external physical phenomena, such as fluid dynamics and heat transfer, and their geometric variations strongly and directly influence performance metrics. Structural and topological optimization problems, where geometry-performance relationships may be indirect or less explicit, constitute a distinct class that is not the main focus of this paper, although selected examples of these methods are briefly referenced in the concluding section.}

The aim \change{here} is to delineate the roles \change{design-space dimensionality reduction} techniques \change{for external shape optimization} play in different engineering applications, emphasizing how they can be harnessed to streamline computational processes, enhance accuracy, and expedite the exploration of optimal designs. Through detailed discussions and case studies, this paper seeks to equip practitioners and researchers with the knowledge to effectively implement dimensionality reduction techniques \change{for shape optimization} in their respective fields, fostering innovation and efficiency in engineering design \cite{li2022machine,serani2024scoping}.


\section{Shape Optimization Problem Formulation}
\change{This section underlines the formulation of the shape optimization problem discussed in the present survey. Specifically, }
consider a manifold $\calG$, which identifies the original/parent shape, whose coordinates in the 3D-space are represented by $\bfg(\bfxi)\in\mathbb{R}^3$; $\bfxi\in\calG$ are curvilinear coordinates defined on $\calG$. Assume that, for the purpose of shape optimization, $\bfg$ can be transformed to a deformed shape/geometry $\bfg'(\bfxi,\bfu)$ by
\begin{equation}
\bfg'(\bfxi,\bfu)=\bfg(\bfxi)+\bfdelta(\bfxi,\bfu) \qquad \forall \bfxi\in\calG
\end{equation} 
where $\bfdelta(\bfxi,\bfu)\in\mathbb{R}^3$ is the resulting shape modification vector, defined by arbitrary shape parameterization or modification method based on fully-parametric or partially-parametric models \cite{harries2019faster}, and $\bfu\in\calU\subset\mathbb{R}^M$ is the design variable vector. 
\change{Figure \ref{fig:deform} illustrates the general concept of shape modification, showing how an original (or parent) geometry $\bfg(\bfxi)$ defined by curvilinear coordinates $\bfxi$, can be transformed (mapped) into a modified geometry $\bfg'(\bfxi,\bfu)$ through a shape modification vector $\bfdelta(\bfxi,\bfu)$. Specifically, this vector captures the geometric displacement (modification) of each point on the surface (or manifold) due to changes in the design variables $\bfu$. This example highlights the flexibility and explicit nature of the shape optimization procedure, clarifying visually the geometric implications of manipulating the parameter vector $\bfu$.}
\begin{figure}[!t]
    \centering
    \includegraphics[width=0.5\columnwidth]{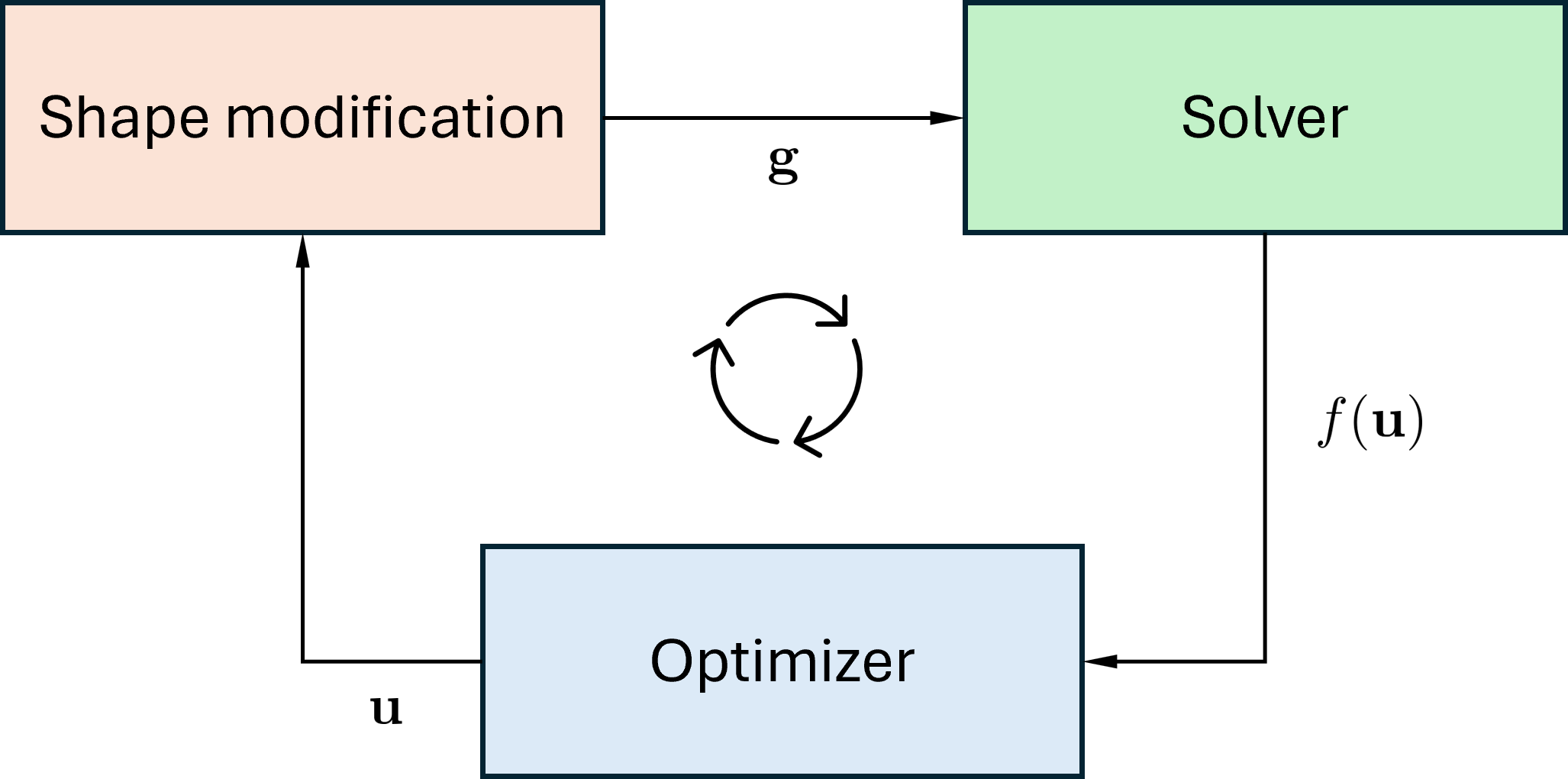}
    \caption{Typical simulation-based design optimization \change{block} diagram for shape optimization}
    \label{fig:xdsm_sbdo}
\end{figure}
The identification of the optimal design $\bfg^\star$ through a general optimization problem reads
\begin{align}\label{eq:soform}
\min_{\bfu} \quad & f(\bfu) \\
\text{subject to} \quad & g_i(\bfu) \leq 0, \quad i = 1, \ldots, I \nonumber \\
\text{and to} \quad & h_j(\bfu) = 0, \quad j = 1, \ldots, J \nonumber\\
\text{and to} \quad & \bfu_l\leq\bfu\leq\bfu_u . \nonumber
\end{align}    
where $f$ is the objective function, $\bfu_l$ and $\bfu_u$ are the lower and upper bounds for the design variable vector, $g_i$ the inequality constraints, and $h_j$ the equality constraints. 
\change{Figure \ref{fig:xdsm_sbdo} provides a high-level flowchart of a typical simulation-based design optimization procedure for shape optimization. The process starts with the initial definition of the design variables and constraints, followed by a loop of design generation (shape modification), numerical simulations (solver), evaluation of the objective and eventual constraint functions, and optimization (optimizer). This iterative cycle is essential to systematically refine the geometry and achieve optimal design solutions that fulfill specified performance criteria.}

\section{Tackling the Curse of Dimensionality}
\begin{figure}[!t]
    \centering
    \includegraphics[width=0.5\columnwidth]{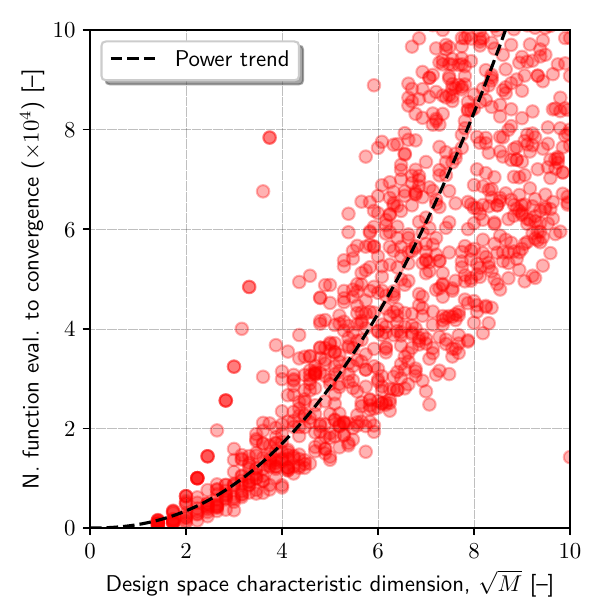}
    \caption{Example of a global optimizer convergence cost conditional to the design space characteristic dimension (data taken from \cite{serani2017random})}
    \label{fig:conv}
\end{figure}

The curse of dimensionality, a term coined by Richard Bellman in 1961, describes the exponential increase in complexity that accompanies adding dimensions to a mathematical \change{problem} \cite{bellman1966dynamic}. In the context of shape optimization, as the dimensionality $M$ of the design variable vector $\mathbf{u}$ \change{(along with, to some extent, the design space bounds)} increases, the volume of the space grows exponentially, making it increasingly difficult to cover it adequately with a finite number of observations. This complexity not only impacts the efficacy of optimization methods, see, e.g., Fig. \ref{fig:conv}, including those relying on surrogate models and multi-fidelity approaches but also complicates surrogate training and uncertainty quantification.

\begin{figure*}[!b]
    \centering
    \includegraphics[width=1\textwidth]{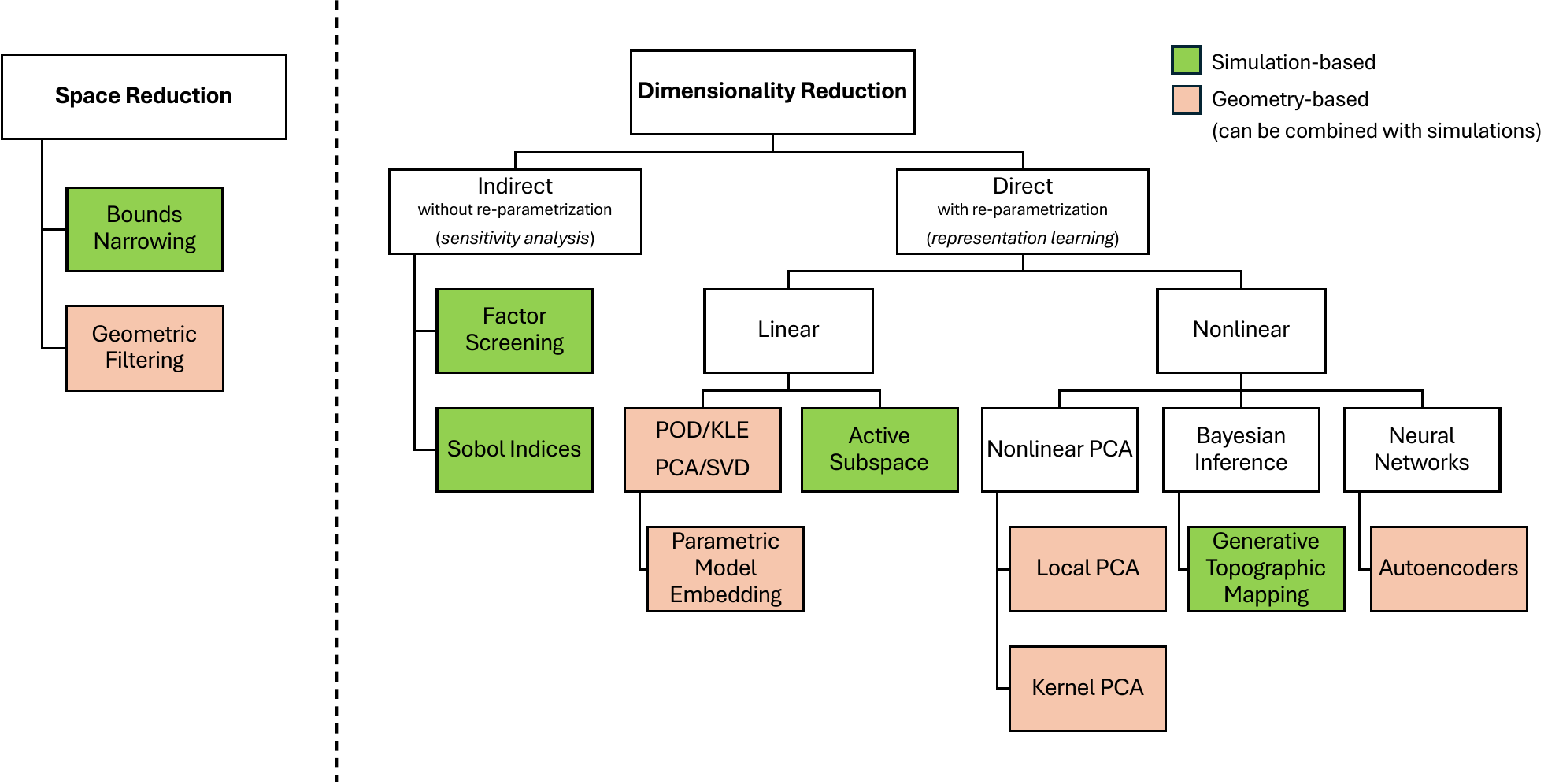}
    \caption{Classification tree for design-space dimensionality reduction in shape optimization}
    \label{fig:class-tree}
\end{figure*}

Figure \ref{fig:class-tree} provides a systematic classification of dimensionality reduction techniques applied in shape optimization. Addressing the curse of dimensionality in such a complex domain requires distinguishing between two principal strategies: space reduction and dimensionality reduction. \change{Each branch clearly highlights the specific methods discussed in this survey, providing readers with a comprehensive framework to navigate the broad landscape of techniques explored in subsequent sections} 

Space reduction techniques focus on reducing the volume of the design space without altering the number of variables involved. This approach typically involves narrowing the bounds within which design variables operate. By constraining the design space to a more manageable volume, these techniques simplify the optimization task without directly reducing the number of design variables. This method is particularly useful in early design phases where broad design constraints are known and can be tightly defined.

In contrast to space reduction, dimensionality reduction techniques decrease the cardinality of the design space—that is, the number of variables themselves. This can be achieved indirectly or directly.

Indirect methods include sensitivity analysis and factor screening techniques, such as one-at-a-time approaches and more complex methods like Sobol indices. These methods identify and eliminate, without re-parameterization, less influential variables based on their impact on output variability. This not only simplifies the number of design variables but also focuses the optimization process on the most significant aspects of the design.

Direct methods re-parameterize the design space by transforming the original high-dimensional space into a lower-dimensional latent space that captures the essential characteristics of the original space. This category includes linear techniques like principal component analysis (PCA) \cite{jolliffe2016principal} and singular value decomposition (SVD), which are discrete representations of proper orthogonal decomposition (POD) \cite{berkooz1993proper}, Karhunen-Loève transform (KLT), or Karhunen-Loève expansion (KLE). Nonlinear approaches extend this concept through methods such as local PCA \cite{kambhatla1997dimension}, kernel PCA \cite{scholkopf1998nonlinear}, Bayesian inference-based techniques like generative topographic mapping (GTM) \cite{bishop1998gtm}, and neural network-based strategies like autoencoders.

Furthermore, design-space dimensionality reduction techniques can be categorized based on their reliance on simulation or geometric data. Simulation-driven techniques, such as indirect methods and some forms of GTM, derive their effectiveness from detailed understanding and manipulation of the objective functions in simulations. Conversely, methods like PCA and its variants primarily utilize geometric data. However, modern approaches often blend these data types, leveraging both geometric and physical simulation data to enhance the applicability and effectiveness of the dimensionality reduction \cite{sobester2013design}.

The choice of dimensionality reduction techniques in shape optimization significantly influences the efficiency and success of the optimization process. By classifying these methods into distinct categories and understanding their underlying principles and applications, engineers and researchers can better select and apply the most appropriate techniques for their specific design challenges. This classification not only aids in comprehending the vast landscape of dimensionality reduction but also highlights the innovative ways these techniques are integrated and applied in the complex field of shape optimization of functional surfaces.

\change{It should be noted that the primary goal of this review is to systematically map which dimensionality reduction methods are currently employed within shape optimization, given the unique severity of the \emph{curse of dimensionality} in this context due to expensive performance evaluations. While certain techniques, particularly those within the PCA family, show specific adaptations tailored for shape optimization problems, others are employed directly, with no significant modifications compared to their general usage in other engineering domains.}

\subsection{Literature overview}
To empirically assess the prevalence and efficacy of design space dimensionality reduction methods, a \change{scoping} review \cite{arksey2005scoping,peters2015guidance,munn2022scoping} was conducted, focusing on their applications across various domains of shape optimization. \change{Scoping reviews, arisen in clinical medicine, use a broad approach for mapping literature and addressing broad research questions without performing articles' quality assessment \cite{sharma2023write}.} 

\change{
\subsubsection{Research Questions} 
The scoping review focuses on the examination of dimensionality reduction methods used in shape optimization, specifically aimed at addressing the challenges associated with high-dimensional design spaces. The following research questions guide this investigation:
\begin{enumerate} 
\item Which dimensionality reduction methodologies are currently applied to manage the curse of dimensionality in shape optimization problems? 
\item What advantages and limitations are associated with applying dimensionality reduction techniques to shape optimization? 
\item What emerging directions and opportunities exist for future research and development in dimensionality reduction for shape optimization?
\end{enumerate}
}

\change{
\subsubsection{Inclusion and Exclusion Criteria} 
The inclusion criteria for this review were carefully defined to identify literature directly addressing shape optimization problems through dimensionality reduction methods. Conversely, the exclusion criteria aimed to ensure the relevance and quality of the selected articles by omitting studies not explicitly related to shape optimization. Additionally, non-peer-reviewed contributions, such as conference papers, abstracts, editorials, and other grey literature, were excluded to maintain academic rigor.
}

\change{
\subsubsection{Databases and Keywords} 
Scopus and Web of Science (WoS) were selected as primary databases for this review due to their extensive coverage of interdisciplinary literature, thus providing a comprehensive and rigorous foundation for identifying relevant studies in shape optimization.
}

\change{
The bibliographic search strategy employed targeted combinations of specific keywords to effectively capture relevant research on dimensionality reduction techniques within shape optimization. The query terms used were: \texttt{(``shape optimi*ation'' OR ``form optimi*ation'' OR ``design optimi*ation'' OR ``geometr* optimi*ation'') AND (``dimension* reduc*'' OR ``space reduc*'' OR ``reduce* dimension*'' OR ``reduced number of dimens*'' OR ``reduced design variable*'' OR ``reduced number of design variable*'')}. These keywords were strategically selected and searched within article titles, abstracts, and keywords sections (using \texttt{TITLE-ABS-KEY} filter for Scopus and \texttt{TS} filter for WoS), ensuring the capture of relevant articles specifically addressing dimensionality reduction methods in shape optimization.
}

\begin{figure}[!b]
    \centering
    \includegraphics[width=0.5\columnwidth]{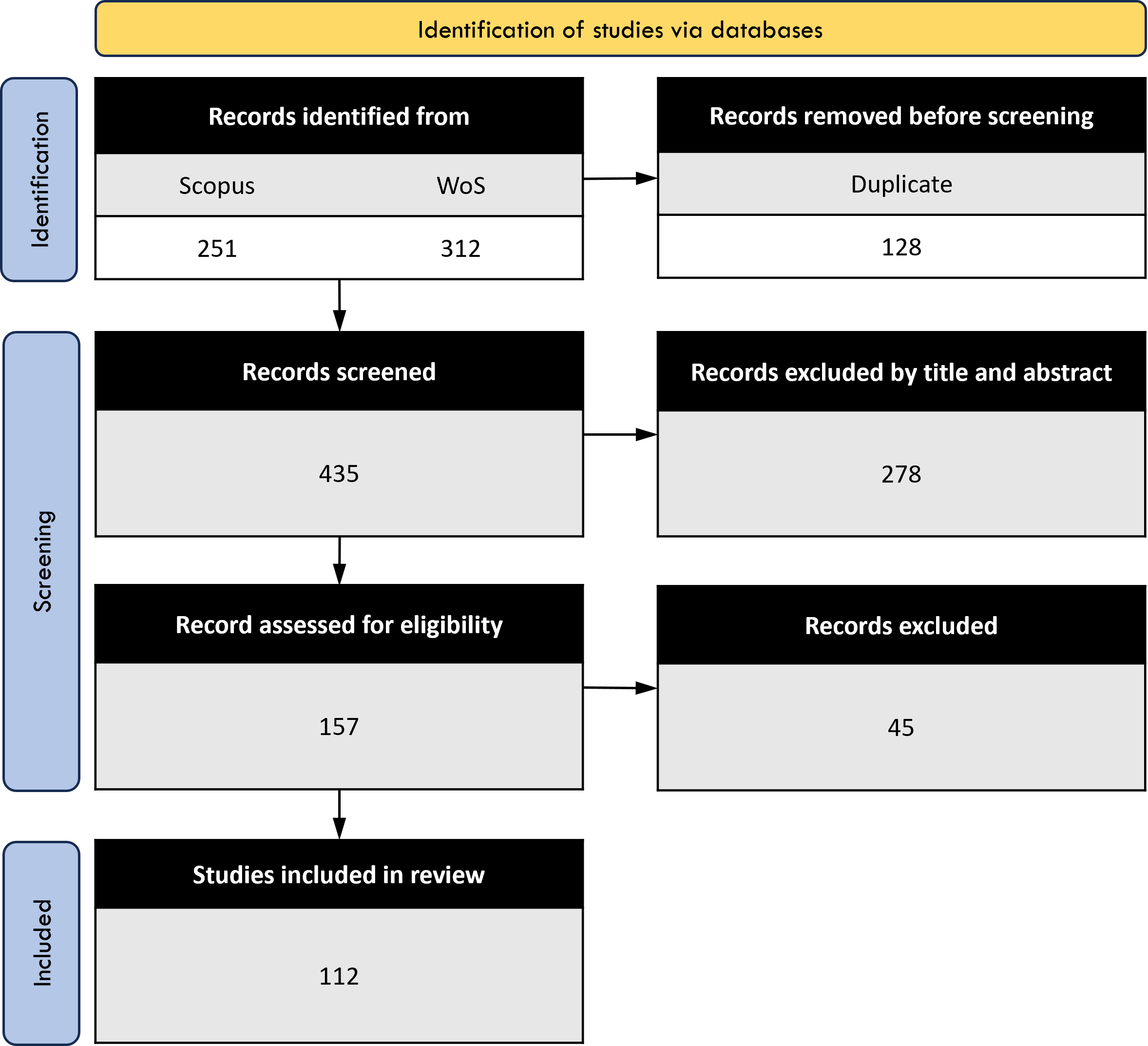}
    \caption{\change{PRISMA flow chart}}
    \label{fig:prisma}
\end{figure}

\change{
\subsubsection{Search Procedure} The search procedure followed the preferred reporting items for systematic reviews and meta-analyses extension for scoping reviews (PRISMA-ScR) guidelines \cite{tricco2018prisma}. The process is depicted in the PRISMA-ScR flow diagram (see Fig. \ref{fig:prisma}), detailing the steps from initial identification to the final selection of articles.
}

\change{
Initially, a total of 435 records were identified through database searches conducted on Scopus and Web of Science (WoS). To further enhance comprehensiveness, reference lists of the included articles were manually examined to identify additional relevant studies.
}

\change{
Following the initial search, a screening phase assessed the relevance of titles and abstracts, resulting in the exclusion of 278 records and retaining 157 potentially relevant articles for detailed review. Subsequently, during the eligibility assessment phase, a full-text analysis was performed, applying specific inclusion and exclusion criteria. Articles were excluded if they did not explicitly address dimensionality or space reduction methods applied to shape optimization problems.
}

\change{
Ultimately, 112 articles were selected, representing a rigorously curated set of literature directly relevant to the focus of this review.
}

\begin{figure*}[!b]
    \centering
    \begin{minipage}[!t]{0.48\textwidth}
    \includegraphics[width=1\columnwidth]{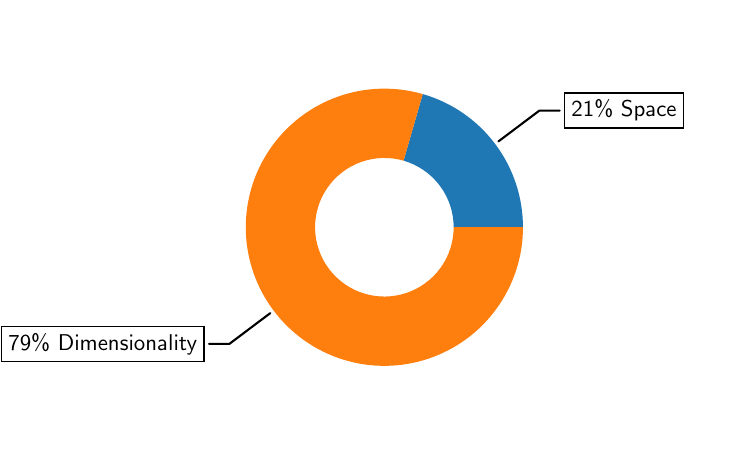}
    \caption{Space versus dimensionality reduction occurrences from shape optimization literature}
    \label{fig:dr_space}
    \end{minipage}
    \hfill
    \begin{minipage}[!t]{0.48\textwidth}
    \includegraphics[width=1\columnwidth]{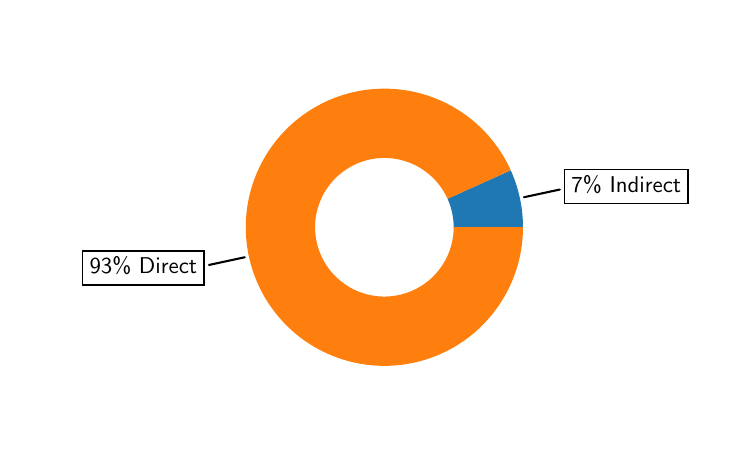}
    \caption{Direct versus indirect dimensionality reduction occurrences from shape optimization literature}
    \label{fig:dr_direct}
    \end{minipage}
\end{figure*}

\change{\subsection{Scoping Review Outcomes}}
The statistical analysis, supported by visual trends from the literature, reveals distinct patterns in the adoption of direct versus indirect methods, the preference for linear over nonlinear approaches, and the allocation between space reduction and dimensionality reduction techniques. This analysis not only reflects the current state-of-the-art but also guides future research directions by identifying areas lacking in-depth exploration and potential methodological synergies.

\begin{figure*}[!b]
    \centering
    \begin{minipage}[!t]{0.48\textwidth}
    \includegraphics[width=1\textwidth]{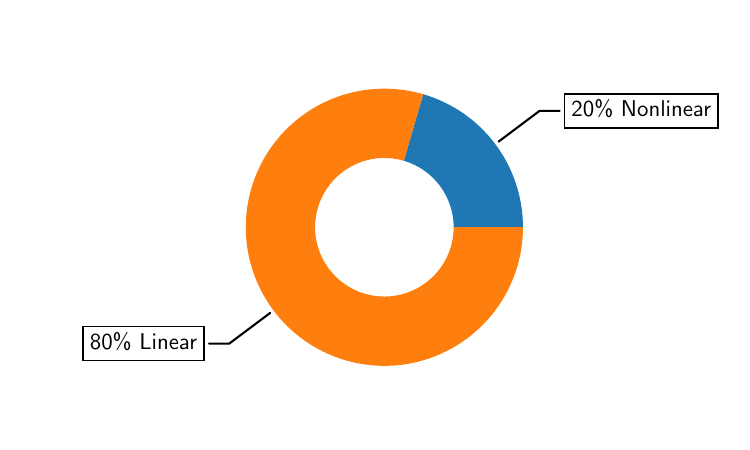}
    \caption{Linear versus nonlinear dimensionality reduction occurrences from shape optimization literature}
    \label{fig:dr_lienar}
    \end{minipage}
    \hfill
    \begin{minipage}[!t]{0.48\textwidth}
    \includegraphics[width=1\textwidth]{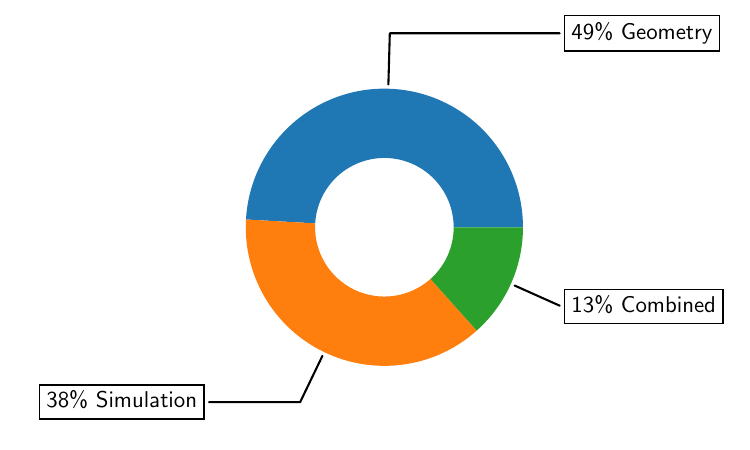}
    \caption{Information used by space and dimensionality reduction methods from shape optimization literature}
    \label{fig:dr_comb}
    \end{minipage}
\end{figure*}

The classification and subsequent statistical overview provide a structured approach to understanding how different dimensionality reduction techniques can be optimally applied to overcome the challenges inherent in high-dimensional design spaces, thus enhancing the overall efficiency of shape optimization processes.

When distinguishing between space reduction and dimensionality reduction, as shown in Fig. \ref{fig:dr_space}, \change{79\% of the studies focus on dimensionality reduction while only 21\% address space reduction}. This highlights a predominant trend towards techniques that reduce the number of variables directly, rather than just condensing the space within which the variables exist.

Considering dimensionality reduction methods, the examination of literature trends underscores a significant leaning towards direct dimensionality reduction methods, with a staggering \change{93\% of reviewed articles opting for direct techniques as opposed to only 7\% utilizing indirect methods}. This suggests a strong preference for methods that actively transform the input data space to achieve dimensionality reduction (see Fig. \ref{fig:dr_direct}).

In terms of methodological choices for direct dimensionality reduction approaches (see Fig. \ref{fig:dr_lienar}), there is an overwhelming preference for linear techniques, with \change{80\% of the cases applying linear methods over 20\% for nonlinear approaches}. This could indicate that the simplicity and computational efficiency of linear methods, such as PCA and SVD, continue to make them popular choices in many engineering applications despite the increasing capability of nonlinear methods to capture more complex data relationships.

Finally, analyzing the approaches based on their reliance on data types, as shown in Fig. \ref{fig:dr_comb}, \change{49\% of the methods are predominantly geometry-based, 38\% are simulation-based, and a smaller segment of 13\% combines both approaches}. This distribution emphasizes the critical role of geometric considerations in the design processes and the significant integration of simulation data to inform these designs.

\begin{figure*}[!t]
    \centering
    \includegraphics[width=1\textwidth]{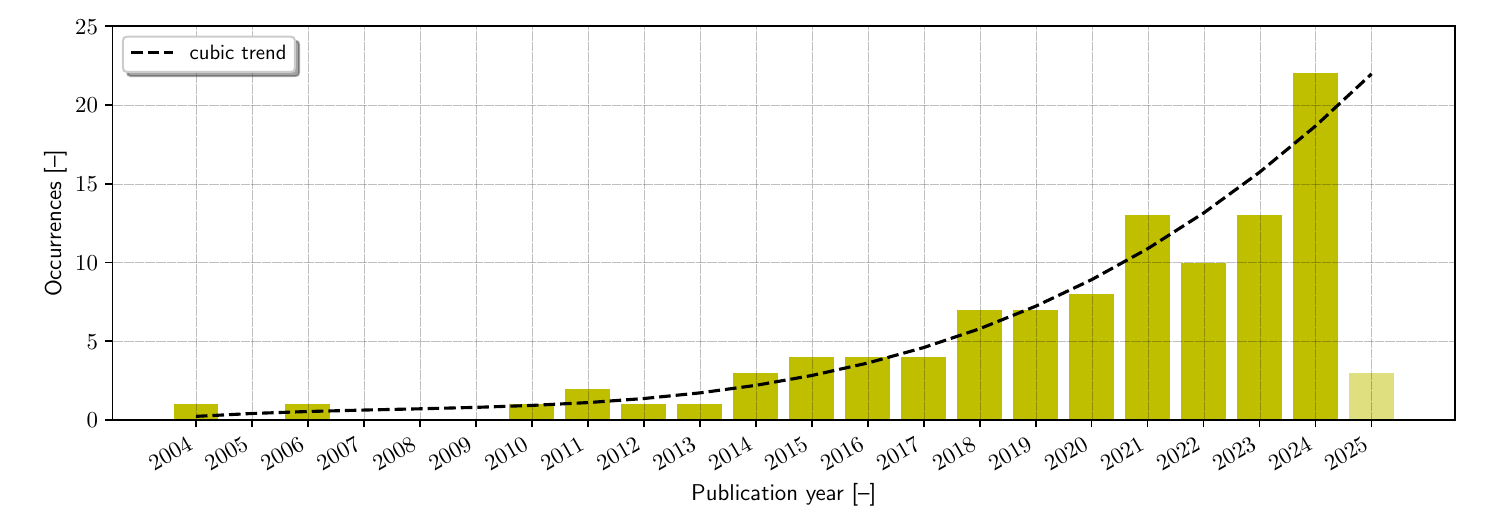}
    \caption{Design-space dimensionality reduction publications trend}
    \label{fig:year}
\end{figure*}

The temporal trend across the years (see Fig. \ref{fig:year}) shows a \change{cubic trend in the increment of} the publications related to dimensionality reduction in shape optimization \change{(note that the data for the year 2025 is partial, as the bibliographic research was conducted on February 26, 2025)}. This growing interest underlines the expanding recognition of the importance of advanced dimensionality reduction techniques in tackling the complexity of modern engineering design problems.

In summary, these trends not only mirror the current preferences and practices in dimensionality reduction but also shed light on the evolving landscape where the integration of advanced computational methods, like machine learning and artificial intelligence, could further revolutionize shape optimization in various industrial applications.

The following sections provide a clear framework and examples for understanding how different space and dimensionality reduction strategies can be effectively utilized in the context of shape optimization, ensuring a more streamlined and focused approach to tackling the intrinsic challenges of high-dimensional optimization problems.

\section{Space Reduction Methods}
\begin{figure*}[!b]
    \centering
    \includegraphics[width=0.75\textwidth]{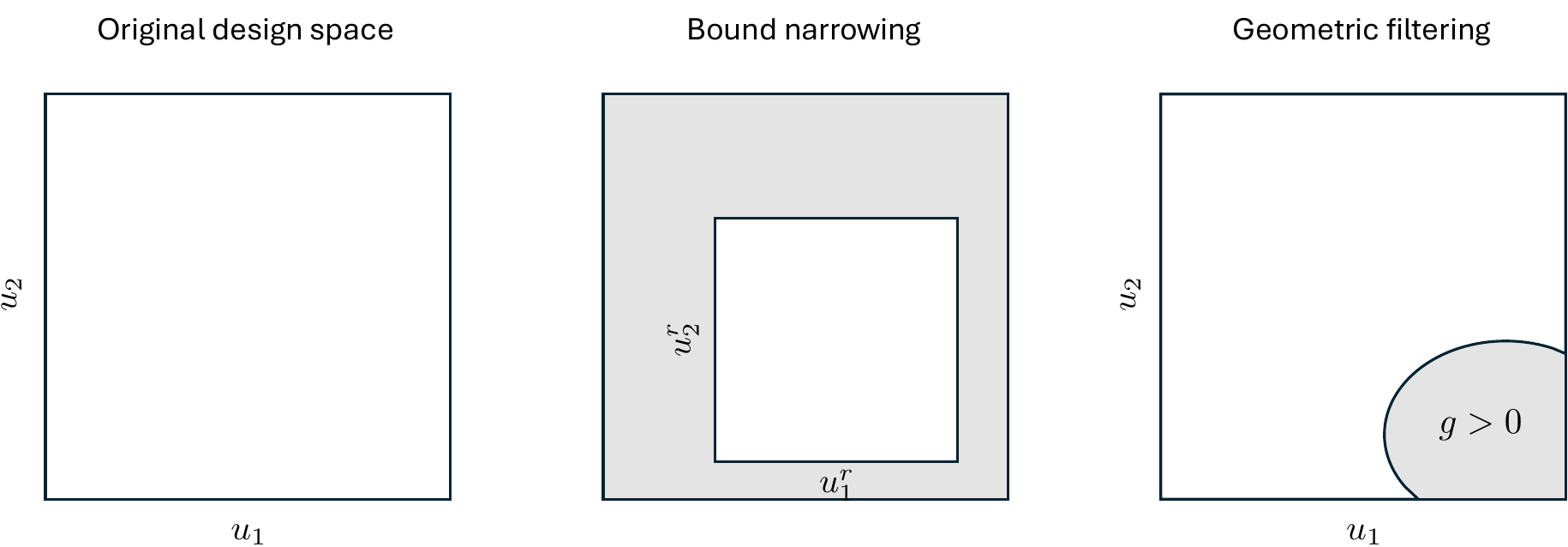}
    \caption{Conceptual representation of space reduction approaches considering $M=2$}
    \label{fig:space_red}
\end{figure*}

Space reduction methods play a crucial role in the field of shape optimization, particularly when dealing with high-dimensional design spaces that challenge traditional optimization techniques. These methods enhance the efficiency and effectiveness of optimization processes by strategically narrowing the focus of the design space. This is achieved either through bounds narrowing, which reduces the breadth of the design space to concentrate on regions most likely to contain optimal solutions, or through geometric filtering, which excludes non-viable or abnormal regions based on sophisticated data-driven models or constraints. 
\change{Figure \ref{fig:space_red} visually contrasts the two main categories of space reduction methods for a simplified two-dimensional design space ($u_1,u_2$).}
Both approaches aim to refine the search area, ensuring that optimization efforts are concentrated on the most promising or feasible parts of the design space. By focusing on these key strategies, the following sections will explore how bounds narrowing and geometric filtering serve as effective tools in simplifying complex optimization tasks, ultimately leading to more innovative and practical engineering solutions.

\subsection{Bounds Narrowing} These methods are instrumental in narrowing the bounds of traditionally expansive design spaces, focusing on smaller, more promising subspaces where optimal solutions are likely to be found.

\change{In aerodynamic shape optimization, bounds narrowing plays a critical role, employing adaptive response surface methodologies and intelligent exploration strategies \cite{koziel2016expedited,koziel2016variable,amrit2017design}. Kriging-based surrogate models, constructed and refined within reduced subspaces, efficiently concentrate on key performance regions, significantly streamlining the optimization process. A variable-fidelity optimization framework, coupled with an intelligent design space reduction strategy, has demonstrated substantial computational time savings while preserving high-fidelity accuracy in airfoil optimization problems \cite{zahir2013variable}. Similarly, a variable-accuracy global metamodel-based optimization approach has been proposed, effectively narrowing the bounds of design experiments and uncertainty quantification for fluid-structure interaction problems \cite{leotardi2016variable}.}

\change{In multidisciplinary optimization contexts, bounds narrowing methods effectively reduce complexity and computational cost. For instance, in optimizing an air-launched satellite vehicle, genetic algorithms combined with simulated annealing successfully narrowed the search space to promising multidisciplinary parameter regions \cite{rafique2011multidisciplinary}. Similar methodologies have also been applied to wing-weight minimization and aero-structural decomposition problems, significantly enhancing computational efficiency \cite{long2015efficient,long2019efficient}.}

\change{In naval and underwater vehicle design, bounds narrowing has been effectively implemented using set-based and hierarchical space reduction methods. Set-based design approaches have been employed in ship design, maintaining flexibility and adapting to evolving requirements within systematically reduced design spaces \cite{hannapel2014implementation}. Dynamic \cite{zheng2021dynamic} and hierarchical \cite{ye2020shape,qiang2022optimization} space reduction frameworks have further advanced ship hull optimization, utilizing data mining techniques such as partial correlation analysis \cite{zheng2021application,zheng2021application2}, self-organizing maps \cite{qiang2022optimization}, and rough set theory \cite{liu2023sequential} to considerably improve optimization efficiency and accuracy. Innovative hierarchical and ensemble surrogate-based optimization methods have significantly improved optimization performance in hydrodynamic shape optimization of blended-wing-body underwater gliders \cite{ye2020shape,ye2024optimization}. Furthermore, bounds narrowing has effectively guided optimization efforts toward critical design regions, demonstrating notable hydrodynamic performance gains for underwater gliders \cite{ran2018two}.}

\change{In structural optimization, a novel method integrating particle swarm optimization with chaotic logistic maps and comprehensive learning strategies has refined dynamically reduced spaces, effectively enhancing the ability to escape local optima and improving the search process for practical-scale engineering problems \cite{van2023chaotic}.}

These applications illustrate the vital role of bounds narrowing in simplifying the optimization process across various domains. By focusing on viable solutions within a confined design space, these strategies ensure more efficient and effective exploration and optimization, underscoring their relevant role in modern engineering design processes.

\subsection{Geometric filtering}

Geometric filtering is a method used to enhance optimization efficiency by excluding abnormal or nonsensical regions within a high-dimensional geometric design space. This approach limits the design space by defining a constraint function to evaluate the abnormality of samples. Various examples highlight how geometric filtering constraints are implemented differently from traditional geometric constraints like distance, thickness, volume, or curvature, as they typically utilize data-based models rather than known equations.

Several case studies where geometric filtering has been applied have been found in the literature. 
An engineering-knowledge-based filtering model using support vector regression to enhance the optimization of a two-dimensional intake duct by incorporating several geometric and flow parameters was developed in \cite{li2012physics}.
An RBF network to eliminate abnormal computer-aided design (CAD) models in nacelle design is employed in \cite{sobester2006supervised}.
An inverse distance weighting to define geometric filtering constraints in modal parameterization of airfoils, which helped train more accurate data-based airfoil analysis models is used in \cite{li2019data}.
Using deep learning, a validity model is created in \cite{li2020efficient} to detect abnormal geometric shapes in airfoils and wing sections, which was trained on a large dataset of labeled airfoils and used in conjunction with generative adversarial networks (GANs) to improve the realism of synthetic airfoils.
The effectiveness of geometric filtering in maintaining innovative designs that enhance aerodynamic performance is also noted, with findings supporting its utility in optimizing designs within transonic regimes \cite{li2021deep}.

Compared to bound narrowing, geometric filtering specifically targets the exclusion of non-viable or abnormal regions from the design space. This method refines the search area by applying advanced data-driven models to define constraints that prevent the exploration of unfeasible or ineffective design solutions. Geometric filtering, therefore, represents a specialized application of space reduction where the reduction is achieved through sophisticated, data-based constraint functions rather than simple boundary adjustments. This approach contributes to more efficient and targeted optimization processes by ensuring that only viable design spaces are considered, thus optimizing computational resources and enhancing the overall effectiveness of the design optimization process.

\section{Indirect Dimensionality Reduction Methods}
\begin{figure*}[!t]
    \centering
    \includegraphics[width=1\textwidth]{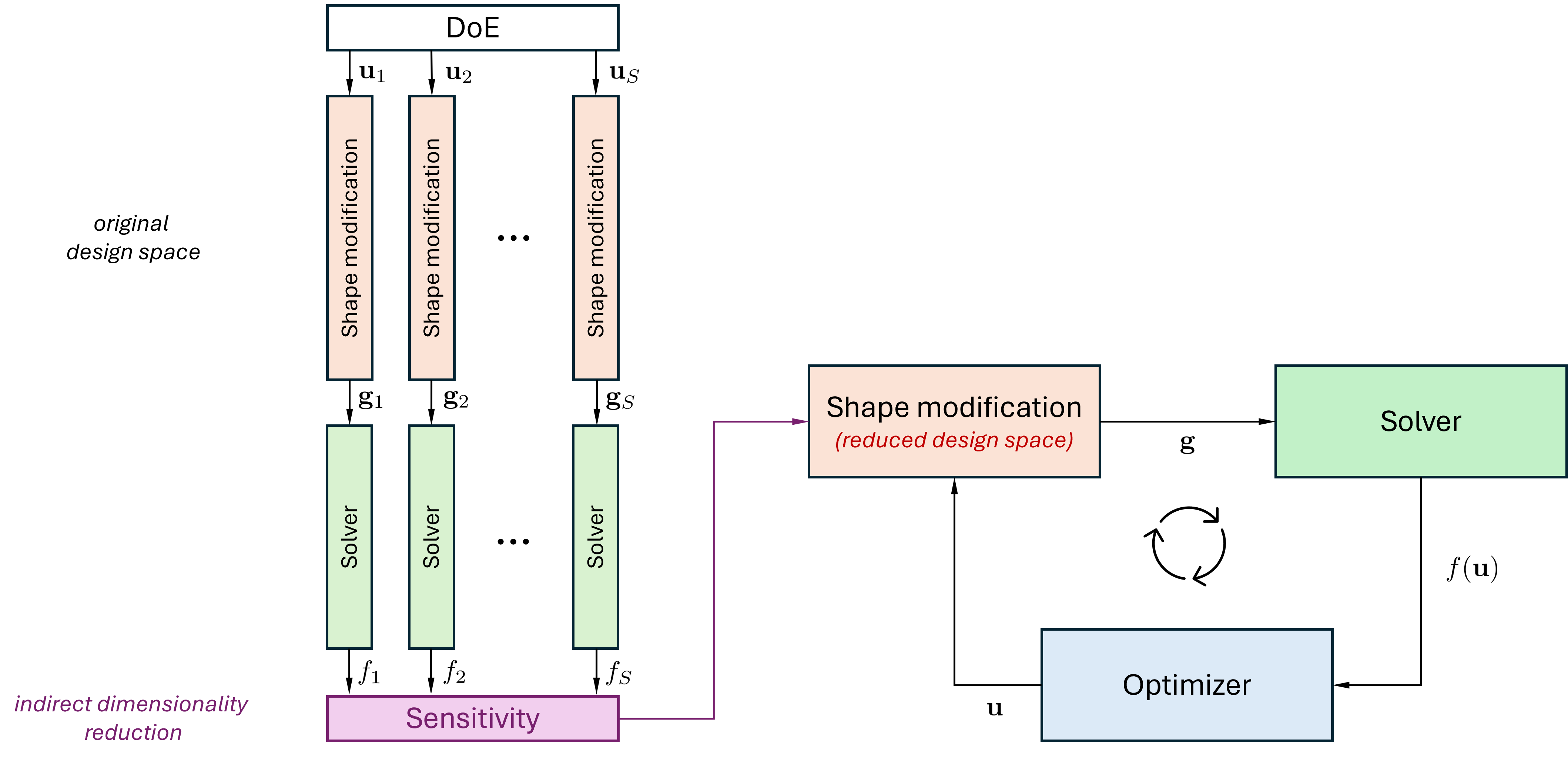}
    \caption{\change{Block diagram for the solution of shape optimization problems through sensitivity methods}}
    \label{fig:block_sensi}
\end{figure*}

Indirect dimensionality reduction methods do not alter the data space itself but focus on identifying the most influential variables and constraints \change{(see Fig. \ref{fig:block_sensi})}. 

Factor screening \cite{montgomery1979factor}, also known as feature selection, helps to identify the most influential factors among a large set of potential variables. By identifying these influential factors early in the optimization process, resources can be focused on refining the optimization of those key factors, potentially leading to faster and more efficient optimization results.

\begin{figure*}[!b]
    \centering
    \includegraphics[width=1\textwidth]{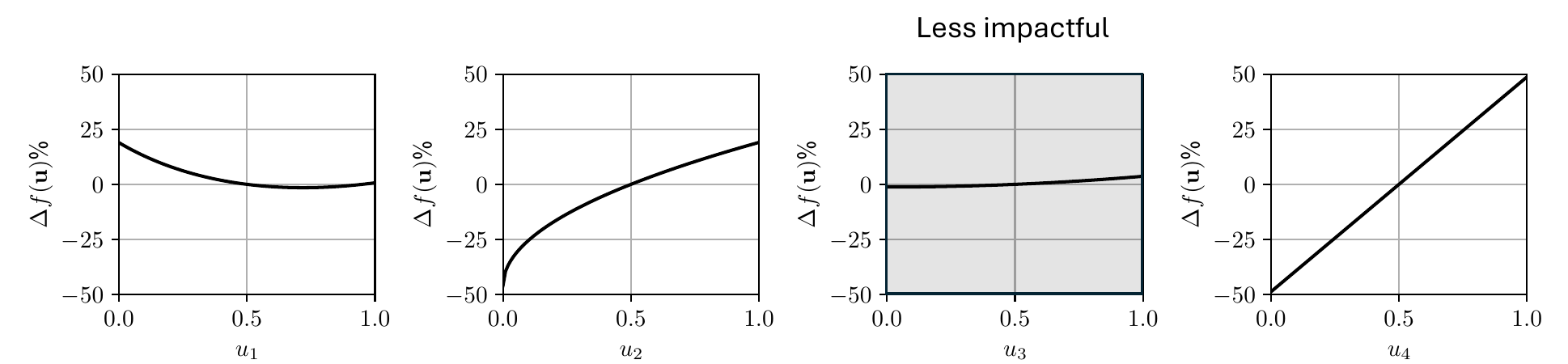}
    \caption{Conceptual representation of sensitivity analysis based on OAT approach for $M=4$; third variable is relatively insensible compared to the other and thus can be removed}
    \label{fig:sensi}
\end{figure*}

In shape optimization, factors may include parameters related to geometry, material properties, boundary conditions, etc. Historically, the simplest factor screening techniques, such as factorial design of experiment and sensitivity analysis, begin with the task of identifying which factors/variables have the most significant impact on the performance or characteristics of the optimized shape. A pivotal step in this preliminary screening process can involve the use of the Morris method for sensitivity analysis \cite{xia2018model}, also known as one-at-a-time (OAT). This approach efficiently sifts through numerous factors to highlight those with the greatest influence on the outcome, by evaluating the change in output resulting from incremental adjustments to each input parameter. This method serves as a practical initial filter, allowing for the identification and retention of key influential factors, while the less impactful ones can be set to a constant value during the subsequent optimization process (see Fig. \ref{fig:sensi}).

Once the influential factors are identified through this combined screening approach, incorporating both the straightforward yet powerful Morris method and other techniques, further optimization efforts can be concentrated on adjusting these key factors to achieve the desired performance goals or characteristics of the optimized shape. This iterative process helps to streamline the optimization process and improve the overall efficiency of shape optimization studies. \change{For instance, the influence of design variables on the objective function is calculated by regression models based on the radial basis functions in \cite{zhang2016sensitivity} and Kriging \cite{xu2024aerodynamic}.}

However, this approach does not always provide the best solution, as it does not fully evaluate the importance that the fixed variables could have during the optimization process, especially when combined with other variables. This limitation can be assessed by variance-based methods like Sobol sensitivity analysis \cite{sobol2001global,gogu2009dimensionality}. The preliminary design of a bulk carrier to enhance its robustness and reliability by focusing only on the most influential random variables, thereby reducing the complexity of the optimization problem, was performed in \cite{wei2019sensitivity}, where polynomial chaos expansion was employed in combination with sensitive analysis to calculate Sobol's indices, by deriving them directly from the coefficients of the polynomial representing the output variable. The approximate local solution of Sobol’s indices was applied to the reliability-based robust design optimization of composite structures in \cite{das2021dimensional}. \change{An improved variance-based sensitivity analysis method was developed specifically for ship hull optimization, where the Sobol method was enhanced by introducing new integrals characterized by reduced variance and employing quasi-random numbers with lower star discrepancy \cite{liu2017improvement}. This enhancement significantly improved the accuracy of sensitivity indices estimation, facilitating an effective reduction of the dimensionality of the hull-form optimization model.}
Nevertheless, the Sobol method operates within a probabilistic framework and consequently requires a statistically significant number of design-space samples. Furthermore, the number of indices increases with the power of the design-space dimensionality, making Sobol indices generally computationally expensive. \change{It is worth noting that surrogate-model-based Sobol analyses can significantly mitigate this computational burden, offering more efficient estimates of sensitivity indices; however, such surrogate-based strategies were not extensively identified within the reviewed shape optimization literature.}


Hence, the need in industrial design is for such dimensionality reduction methods that can capture, in a reduced-dimensionality space, the underlying most promising directions of the original design space, preserving its relevant features and thereby enabling an efficient and effective optimization in the reduced space. 

The remedy can be found in those direct dimensionality reduction techniques classified as unsupervised learning, feature extraction, or also representation learning \cite{bengio2013representation}, capable of learning relevant hidden structures of the original design-space parameterization.

\section{Direct Dimensionality Reduction Methods: Representation Learning}
Representation learning operates at the intersection of machine learning and optimization, presenting a transformative approach to handling complex data structures efficiently. In shape optimization, 
these techniques, known also as \textit{modal parameterization} \cite{li2022machine}, play a fundamental role in extracting insightful information from high-dimensional shape spaces, fostering more effective exploration and manipulation of shapes.
The primary aim of representation learning in shape optimization is to acquire useful representations or embeddings or re-parameterization of the shape modification vectors $\bfu$ that encapsulate their essential characteristics or features. Through re-parameterization, representation learning allows to reformulate problem in Eq. \ref{eq:soform} as follows 
\begin{align}\label{eq:roform}
\min_{\bfx} \quad & f(\bfx) \\
\text{subject to} \quad & g_i(\bfx) \leq 0, \quad i = 1, \ldots, I \nonumber \\
\text{and to} \quad & h_j(\bfx) = 0, \quad j = 1, \ldots, J \nonumber\\
\text{and to} \quad & \bfx_l\leq\bfx\leq\bfx_u . \nonumber
\end{align}    
where $\bfx\in\mathbb{R}^N$ embeds the main features of the original design-space (with $N<M$), and $\bfx_l$ and $\bfx_u$ are its lower and upper bounds, respectively, leading to faster convergence and improved scalability.

Within the context of dimensionality reduction, representation learning encompasses a diverse array of both linear and nonlinear techniques aimed at learning compact and informative representations of the original shape parameterization.

\subsection{Linear dimensionality reduction methods}
Linear dimensionality reduction methods are commonly employed in shape optimization to identify dominant modes of variation and facilitate dimensionality reduction while preserving significant features. Following subsections describes the main methodologies and their mathematical assumptions, along with the relevant references.

\change{\subsubsection{Principal component analysis, from KLE/POD}}

The problem in Eq. \ref{eq:soform} can be considered as affected by an epistemic uncertainty, where $\bfu$ can be assumed as an uncertain/random parameter, thus the optimal design $\bfg^\star$ exists but, before going through the optimization procedure, is unknown. Accordingly, a probability density function $p(\bfu)$, representing the degree of belief in finding the optimal solution in a certain region of the design space, can be assigned to the design variable vector $\bfu$. This probability density is arbitrary and may be chosen based on the designer's previous knowledge and past experience. If previous knowledge is not available or deemed not useful for the optimization study, a uniform distribution may be used, which gives the same probability to all designs in the given domain to be selected as optimal decisions. Once $p(\bfu)$ is defined, the shape modification vector $\bfdelta$ goes stochastic and can be studied as a random field, by using the KLE, equivalent to POD. 

Consider $\bfdelta(\bfxi,\bfu)$ as part of a Hilbert space $L_\rho^2(\calG)$, characterized by the generalized inner product
\begin{equation}\label{eq:inprod}
(\bfa,\bfb)_\rho=\int_\calG \rho(\bfxi)\bfa(\bfxi)\cdot\bfb(\bfxi)\dxi
\end{equation}
along with the associated norm $\|\bfa\|=(\bfa,\bfa)_\rho^{1/2}$, where $\rho(\bfxi)\in\mathbb{R}$ denotes an arbitrary weight function. Examining all potential realizations of $\bfu$, the corresponding mean vector of $\bfdelta$ is given by
\begin{equation}
\left\langle\bfdelta\right\rangle=\int_\calU \bfdelta(\bfxi,\bfu)p(\bfu)\du
\end{equation}
while the corresponding geometrical variance is expressed as
\begin{equation}\label{eq:sigma2}
\sigma^2=\left\langle\|\hat{\bfdelta}\|^2\right\rangle=\iint\limits_{\calU,\calG} \rho(\bfxi)\hat{\bfdelta}(\bfxi,\bfu)\cdot\hat{\bfdelta}(\bfxi,\bfu)p(\bfu)\dxi\du
\end{equation}
where $\hat{\bfdelta}=\bfdelta-\left\langle\bfdelta\right\rangle$, with $\left\langle\cdot\right\rangle$ representing the ensemble average obtained by integrating $\bfu$ over $\calU$.

The objective of POD/KLE is to identify an optimal basis of orthonormal functions for the linear representation of $\hat{\bfdelta}$:
\begin{equation}\label{eq:delta}
\hat{\bfdelta}(\bfxi,\bfu)\approx\sum_{k=1}^N x_k(\bfu)\bfphi_k(\bfxi)
\end{equation}
where
\begin{equation}\label{eq:xred}
x_k(\bfu)=\left(\hat{\bfdelta},\bfphi_k\right)_\rho=\int_\calG \rho(\bfxi)\hat{\bfdelta}(\bfxi,\bfu)\cdot\bfphi_k(\bfxi)\dxi
\end{equation}
constitute the basis-function components usable as new (reduced) design variables. The optimality condition associated with POD/KLE pertains to the geometric variance retained by the basis functions through Eq. \ref{eq:delta}. Combining Eqs. \ref{eq:sigma2}--\ref{eq:xred} yields
\begin{equation}
\sigma^2=\sum_{k=1}^\infty\left\langle\left(\hat{\bfdelta},\bfphi_k\right)_\rho^2\right\rangle
\end{equation}

The basis retaining the maximum variance consists of those $\bfphi$ solutions to the variational problem
\begin{eqnarray}\label{eq:varp}
\underset{\bfphi\in\calL_\rho^2(\calG)}{\rm maximize} & \calJ(\bfphi)=\left\langle\left(\hat{\bfdelta},\bfphi_k\right)_\rho^2\right\rangle\\
\mathrm{subject , to} & \left(\bfphi,\bfphi\right)_\rho^2=1\nonumber
\end{eqnarray}
as detailed in \cite{diez2015design}
\begin{align}\label{eq:bie}
\calL\bfphi(\bfxi)&=\int_\calG \rho(\bfxi')\left\langle\hat{\bfdelta}(\bfxi,\bfu)\otimes\hat{\bfdelta}(\bfxi',\bfu)\right\rangle\bfphi(\bfxi')\dxi' \nonumber\\
&=\lambda\bfphi(\bfxi)
\end{align}
where $\calL$ represents the self-adjoint integral operator whose eigensolutions define the optimal basis functions for the linear representation of Eq. \ref{eq:delta}. Consequently, its eigenfunctions (modes) $\{\bfphi_k\}_{k=1}^\infty$ are orthonormal and form a complete basis for $\calL_\rho^2(\calG)$. Additionally, it may be proven that
\begin{equation}
\sigma^2=\sum_{k=1}^\infty\lambda_k \quad \mathrm{with} \quad \lambda_k=\left\langle x_k^2\right\rangle
\end{equation}
where the eigenvalues $\lambda_k$ represent the variance retained by the associated basis function $\bfphi_k$, through its component $x_k$. Finally, the solution $\{\bfphi_k\}_{k=1}^\infty$ of Eq. \ref{eq:varp} are used to construct a reduced dimensionality representation of the original design space; defining the desired confidence level $l$, with $0<l\leq1$, the number of reduced design variables $N$ in Eq. \ref{eq:delta} is chosen such that
\begin{equation}\label{eq:lambda}
\sum_{k=1}^N \lambda_k\geq l\sum_{k=1}^\infty \lambda_k=l\sigma^2 \quad \mathrm{with} \quad \lambda_k\geq\lambda_{k+1}
\end{equation}

\begin{figure*}[!b]
    \centering
    \includegraphics[width=1\textwidth]{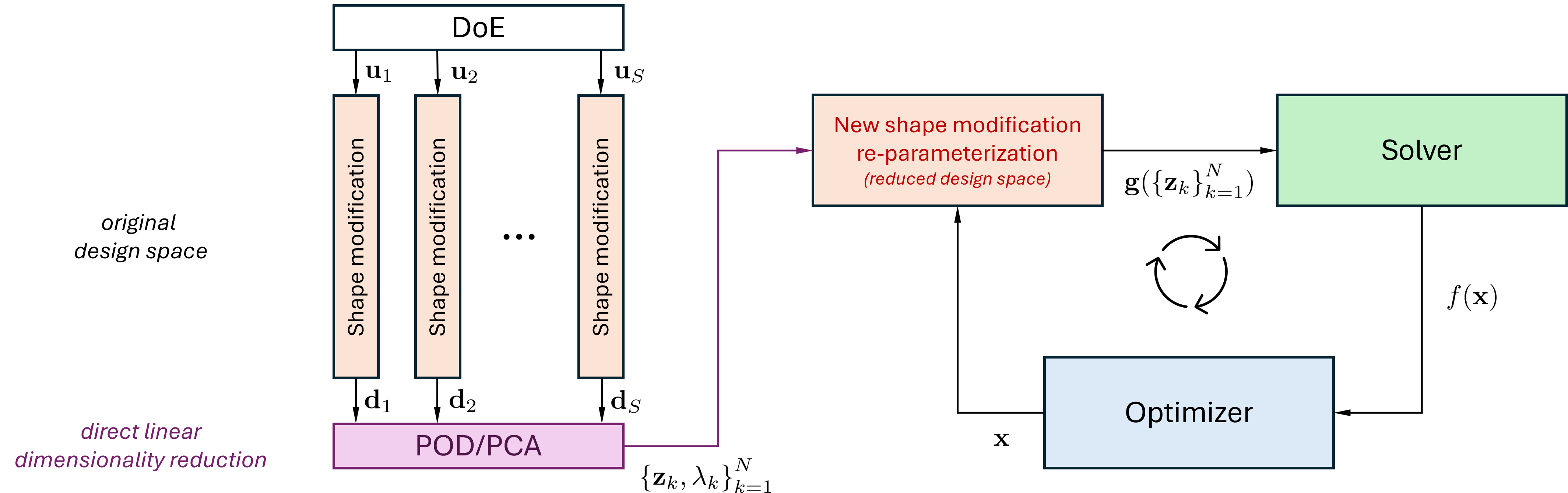}
    \caption{\change{Block diagram for the solution of shape optimization problems through PCA}}
    \label{fig:xdsm_pca}
\end{figure*}
%
%
\begin{figure*}[!t]
    \centering
    \includegraphics[width=0.59\textwidth]{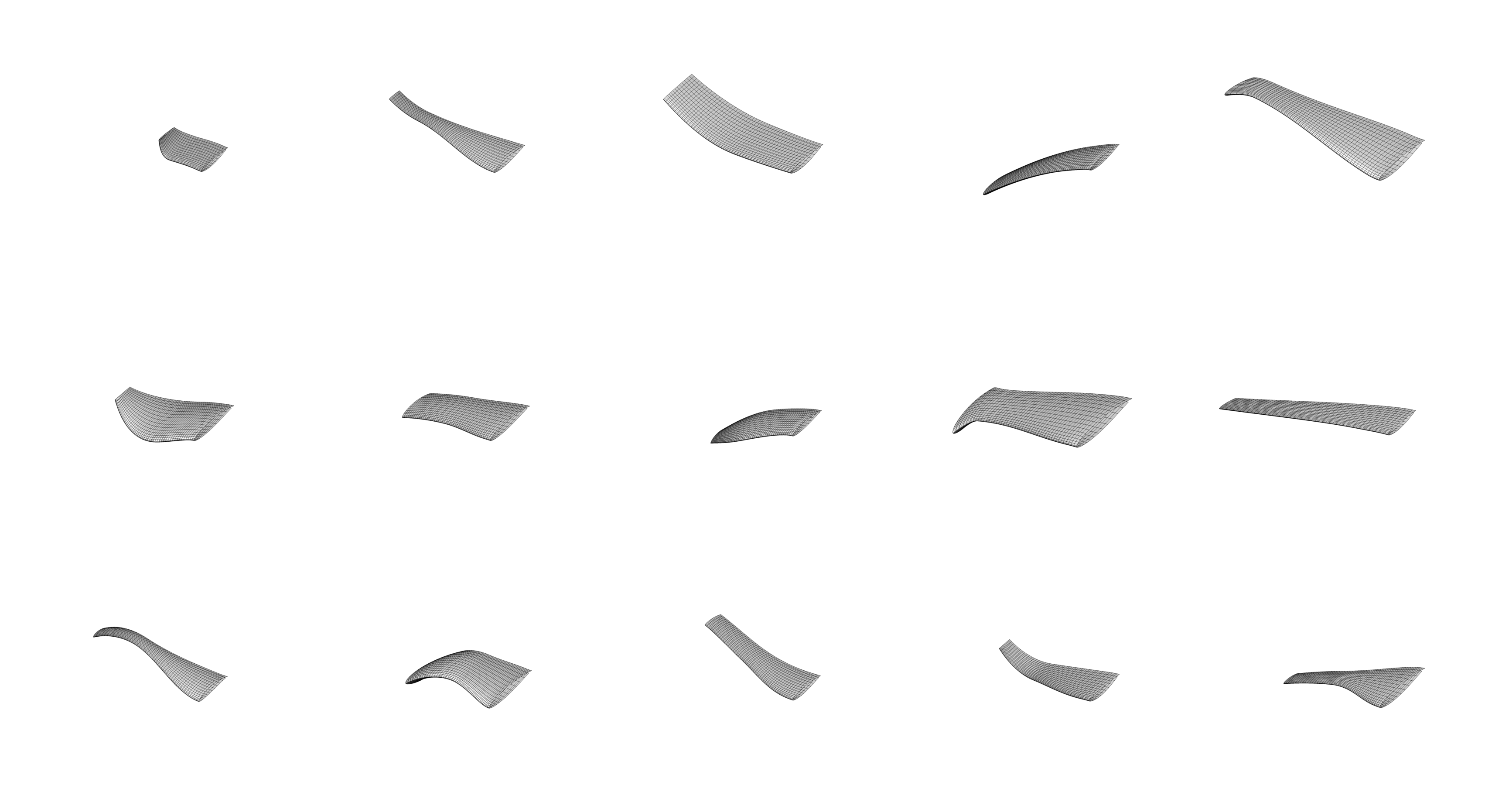}
    \includegraphics[width=0.4\textwidth]{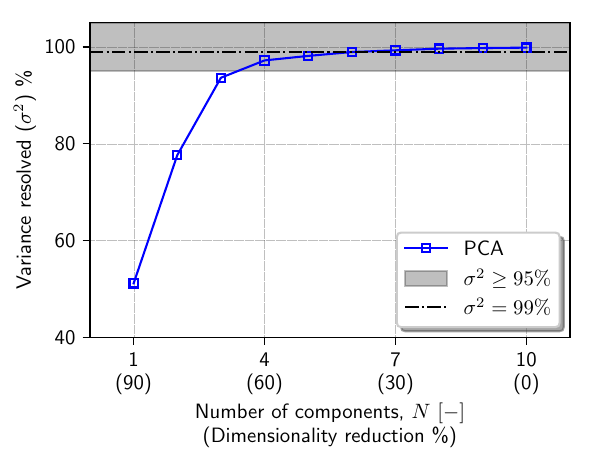}
    \caption{Example of (left) random samples of the design space for the optimization of a hydrofoil \cite{solak2023hydrofoil} and (right) dimensionality reduction results in terms of cumulative geometric variance retained by the PCA components}
    \label{fig:pca}
\end{figure*}
Discretizing the domain $\calG$ into $L$ elements, each with a measure denoted as $\Delta\calG_i$ (where $i=1,\dots,L$), and sampling the input space $\calU$ using a statistically convergent number of realizations, denoted as $S$, results in a set $\{\bfu_k\}_{k=1}^S$ drawn from the distribution $p(\bfu)$. These realizations are organized into a data matrix $\mathbf{D}$ of dimensionality $\left[3L\times S\right]$, given by:

\begin{equation}\label{eq:data}
\mathbf{D}=\left[
\begin{array}{ccc}
\hat{d}_{1,\xi_1}(\bfu_1) & & \hat{d}_{1,\xi_1}(\bfu_S)\\
\vdots & & \vdots\\
\hat{d}_{L,\xi_1}(\bfu_1) & & \hat{d}_{L,\xi_1}(\bfu_S)\\
\hat{d}_{1,\xi_2}(\bfu_1) & & \hat{d}_{1,\xi_2}(\bfu_S)\\
\vdots & \dots & \vdots\\
\hat{d}_{L,\xi_2}(\bfu_1) & & \hat{d}_{L,\xi_2}(\bfu_S)\\
\hat{d}_{1,\xi_3}(\bfu_1) & & \hat{d}_{1,\xi_3}(\bfu_S)\\
\vdots & & \vdots\\
\hat{d}_{L,\xi_3}(\bfu_1) & & \hat{d}_{L,\xi_3}(\bfu_S)\\
\end{array}
\right]
\end{equation}

Here, $\hat{d}_{i,\xi_j}$ denotes the deviation from the mean geometry of $j$-th component of the discretized shape modification vector $\hat{\mathbf{d}}_i(\bfu_k)$ associated with element $i$ of geometry $k$. By considering this discretization and employing spectral analysis, the integral problem of Eq. \ref{eq:varp} is transformed into the generalized PCA of the data matrix $\mathbf{D}$:

\begin{equation}\label{eq:pca}
\mathbf{A {G}WZ}=\mathbf{Z}\boldsymbol{\Lambda} \qquad \mathrm{with} \qquad \mathbf{A}=\frac{1}{S}\mathbf{DD}^\mathsf{T}
\end{equation}

Here, $\mathbf{Z}$ and $\boldsymbol{\Lambda}$ represent the eigenvectors and eigenvalues matrices of $\mathbf{A{G}W}$. The matrix $\mathbf{G}$, defined as a block diagonal matrix with dimensionality $\left[3L\times 3L\right]$, and $\mathbf{W}$, similarly a block diagonal matrix of dimensionality $\left[3L\times 3L\right]$, play crucial roles in this transformation. Specifically, $\mathbf{G}$ consists of block matrices $\mathbf{G}_k$ (for $k=1,2,3$), each being a diagonal matrix:

\begin{equation}
\mathbf{G}_k=\mathrm{diag}\left(\Delta\calG_1, \dots, \Delta\calG_L\right)
\end{equation}

Here, $\Delta\calG_i$ denotes the measure of the $i$-th element. Similarly, $\mathbf{W}$ is a block diagonal matrix composed of diagonal matrices $\mathbf{W}_k$:

\begin{equation}\label{eq:weights}
\mathbf{W}_k=\mathrm{diag}\left({\rho_1}, \dots, {\rho_L}\right)
\end{equation}

Here, $\rho_i$ (for $i = 1,\dots, L$) denotes the arbitrary weight assigned to each element. It's noteworthy that if all elements of $\calG$ have equal measure $\Delta\calG$ and are assigned the same weight $\rho$, then the problem reduces to the standard PCA of $\mathbf{D}$, 
\begin{equation}
\label{eq:pca_standard}
\mathbf{AZ}=\mathbf{Z}\boldsymbol{\Omega} \qquad \mathrm{with} \qquad \mathbf{\Omega}=\frac{1}{\Delta\calG \, \rho}\mathbf{\Lambda}
\end{equation}

Moreover, the discretization of $\calG$ can be arbitrary and doesn't need to align with the discretization of the solver used for shape optimization. However, it can consider the physics of the problem by assigning arbitrary weights to the elements where significant physical phenomena occur.

Lastly, it's worth mentioning that both problems \ref{eq:pca} and \ref{eq:pca_standard} result in the matrix $\mathbf{Z}$ containing the discrete representation $\bfz_k$ of the desired eigenfunctions $\bfphi_k$, that can be used as the basis for the reconstruction of shape modification vector $\hat{\mathbf{d}}$ from the reduced design variables $\bfx$ as
\begin{equation}
\hat{\mathbf{d}}(\bfu)\approx \sum_{k=1}^N x_k(\bfu)\bfz_k.
\end{equation}
The \change{block} diagram shown in Fig. \ref{fig:xdsm_sbdo} can be consequently transformed as shown in Fig. \ref{fig:xdsm_pca}, \change{where it is shown how the dimensionality reduction procedure is performed upfront the optimization loop and fed it with the new parameterization based on PCA eigenvectors and eigenvalues.}

POD/PCA have been crucial in enhancing the design of airfoils and wings, streamlining external aerodynamic optimizations for both subsonic \cite{yonekura2014shape,cinquegrana2018investigation,li2019data}, transonic \cite{toal2010geometric, poole2015metric, poole2017high, masters2017geometric, allen2018wing, yasong2018global, li2019data, wu2019benchmark, jun2020application, li2021adjoint, yamazaki2020efficient, poole2022efficient, wu2024high}, and hypersonic \cite{li2024aerodynamic} regimes of airfoils and wings, as well as for internal aerodynamics of compressors \cite{zhang2018multidisciplinary, yanhui2019performance}, turbines \cite{hu2022dimension}, and nozzles \cite{yang2022optimal}. These methods have facilitated efficient parameterization and reduced drag, optimizing aerodynamic performance for aircraft \change{and watercraft, including blended wing body \cite{yang2024aerodynamic} and} underwater autonomous \cite{li2022low} vehicles.
In marine applications, KLE/PCA has improved hydrodynamic performance by optimizing the designs of mono- \cite{harries2019faster,ballarin2019pod,d2020design,chang2023research,zhang2024geometric} and multi-hulls \cite{diez2015design, chen2015high,harries2021application}, hull vanes \cite{ccelik2021reduced}, propellers \cite{gaggero2020reduced}, \change{and turbines \cite{masood2021machine}}. These applications have demonstrated effective reductions in design variables and computational costs while improving resistance and efficiency. Naval hull optimizations \cite{padula2024generative} have particularly benefited from generative models and advanced PCA techniques \cite{d2024generative} that manage design deformations and enhance performance. Figure \ref{fig:pca} shows an example of random samples (used for the construction of the data matrix $\mathbf{D}$) of $M=10$ parametric design space for the optimization of a hydrofoil and the associated dimensionality reduction result \cite{solak2023hydrofoil} \change{in terms of geometric variance resolved as a function of the number of reduced design variables}.

These dimensionality reduction strategies have not only optimized computational resources but also enhanced the quality and efficiency of designs across aerodynamics and marine engineering, proving essential for modern engineering optimization processes.

\begin{figure*}[!t]
    \centering
    \includegraphics[width=1\textwidth]{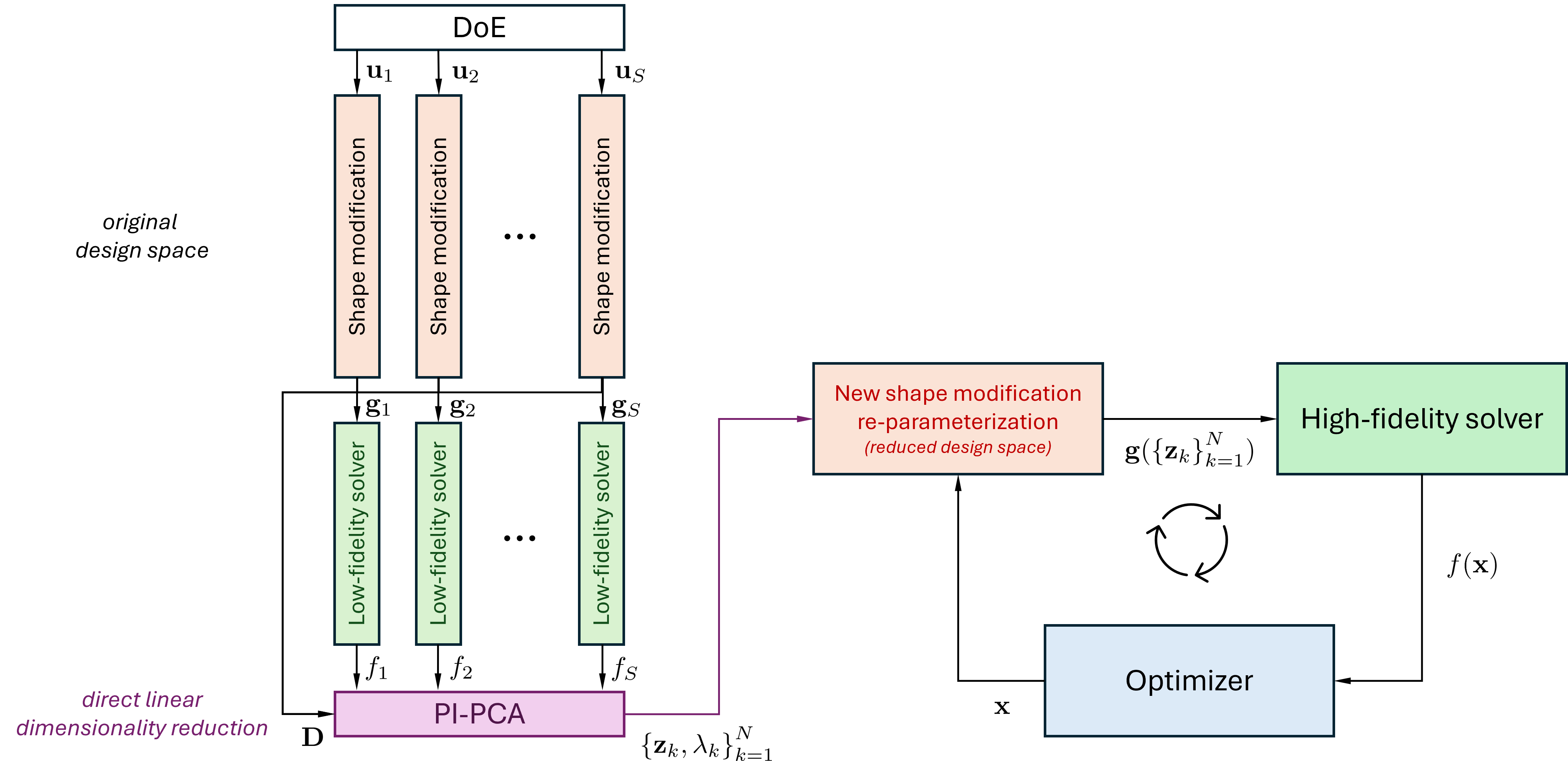}
    \caption{\change{Block diagram for the solution of shape optimization problems through physics-informed PCA}}
    \label{fig:block_pipca}
\end{figure*}
\subsubsection{Physics-informed principal component analysis}
Design-space dimensionality reduction using PCA traditionally relies solely on geometric information. However, by augmenting the data matrix $\mathbf{D}$ with physical data derived from numerical simulations, the method's applicability and effectiveness can be substantially enhance. This augmentation can include both distributed values, such as velocity or pressure fields, and lumped parameters like, e.g., resistance or efficiency. By integrating these additional data types, the PCA not only captures geometric variability but also reflects critical physical phenomena that could influence the design's performance. \change{Figure \ref{fig:block_pipca} shows how PCA-based workflow (see Fig. \ref{fig:xdsm_pca}) can be augmented by computing low-fidelity physical quantities for each design variant of the original design space upfront the optimization loop}.

Incorporating physics-informed parameters allows for a more holistic understanding of the design space by bridging the gap between purely geometric insights and their functional outcomes. 
For example, PCA integrated with physical data has been demonstrated to effectively map the design space of waterjet propulsion systems, considering shape and flow characteristics essential for optimizing performance \cite{liu2021linear}. Similarly, aerodynamic shape optimization of turbomachinery has benefited from the integration of aerodynamic sensitivities into PCA, improving uncertainty quantification by capturing geometry variations critical to performance metrics \cite{li2023new,mingzhi2024breaking} and enhancing the efficiency of deviation modeling in compressors \cite{mingzhi2024breaking}.
Furthermore, a metric-based POD parameterization method has demonstrated significant improvements in aerodynamic shape optimization efficiency, employing data-based filtering strategies to refine design spaces specifically towards better aerodynamic characteristics \cite{zhang2023efficient}. In rotorcraft applications, aerodynamic and aeroacoustic optimization has been substantially advanced by integrating POD parameterization with physical aerodynamic data, enabling accurate and efficient exploration of design spaces sensitive to performance and noise constraints \cite{liu2024aerodynamic}.
In naval hydrodynamics, PCA combined with physics-informed simulation data has been successfully employed to analyze ship hydrodynamics, assessing how shape variations influence wave resistance and other critical performance parameters \cite{serani2020assessing}. Moreover, recent studies explored a combined shape- and physics-based PCA formulation that significantly improves hydrodynamic shape optimization efficiency by integrating low-fidelity solver outputs into preliminary design phases, thereby efficiently identifying key design modifications contributing to enhanced deterministic \cite{diez2023design} and stochastic performance \cite{serani2022hull} under multi-fidelity conditions.

Additionally, \cite{khan2022shape, khan2022geometric, masood2023shape, kostas2023machine,masood2024generative} introduced geometric moments as cost-effective, physics-related descriptors that substitute direct simulation data, further enriching the PCA model without the computational cost of numerical simulations. This adaptation enables the efficient exploration of design spaces by focusing on features that are both geometrically significant and closely related to performance metrics.

These enhancements to PCA, \change{leveraging physics-informed methods}, not only increase the accuracy of the dimensionality reduction but also enable designers to make more informed decisions early in the design process, potentially reducing the number of high-fidelity simulations required and thereby lowering computational costs.

\subsubsection{Parametric model embedding}
If the PCA-based dimensionality reduction procedure exclusively employs information on the shape modification vector, it does not inherently offer a means to revert to the initial design variables from the latent space (i.e., the reduced dimensionality representation of the original shape parameterization). This situation presents two significant challenges: 
\begin{enumerate}
    \item relying solely on shape-modification-based PCA necessitates the development of a new shape modification method based on the PCA eigenvectors; 
    \item depending on the bounds applied to the reduced design variables, there is no assurance that the shape generated using PCA eigenvectors indeed resides within the original design space, potentially leading to design infeasibilities.
\end{enumerate}
A systematic approach is imperative to establish an explicit relationship between original and reduced design variables, a fundamental aspect for the widespread utilization of PCA-based design-space dimensionality reduction methods in industrial contexts. In such settings, design parameters frequently conform to well-established parametric models (e.g., CAD models), which serve as essential inputs for numerical solvers (e.g., for isogeometric analysis). Addressing this challenge entails extending the methodology to address the \emph{pre-image problem} \cite{gaudrie2020modeling}, where the objective is to construct a direct mapping of the original variables from the reduced design space, effectively creating a reduced-dimensionality embedding of the original design parameterization. Despite recent efforts to tackle the pre-image problem using interpolation/regression methods, such as the back-transformation approach proposed in \cite{bergmann2018massive} and the Gaussian process model in \cite{gaudrie2020modeling}, these approaches typically lack explicit mathematical relationships between the original parameterization and its reduced-dimensionality representation and may pose challenges in application.

A novel approach, termed parametric model embedding (PME), recently introduced in \cite{serani2023parametric}, offers a straightforward solution to embed the original parameters of parametric models in a reduced-dimensional representation. PME facilitates a direct pathway to revert to the original parameterization from the reduced-dimensionality latent space, enabling the super-parameterization of the original parametric model. Unlike classical sensitivity analysis, where the most critical variables are identified from a design of experiments \cite{harries2019faster}, PME's superparameters exist within the reduced-dimensionality space and do not merely constitute a subset of the original variables in the parametric space.

PME extends the PCA formulation by incorporating a generalized feature space encompassing shape modification and design variable vectors, along with a generalized inner product, aimed at addressing prescribed design variability through judicious selection of the latent dimensionality. 

\begin{figure*}[!b]
    \centering
    \includegraphics[width=1\textwidth]{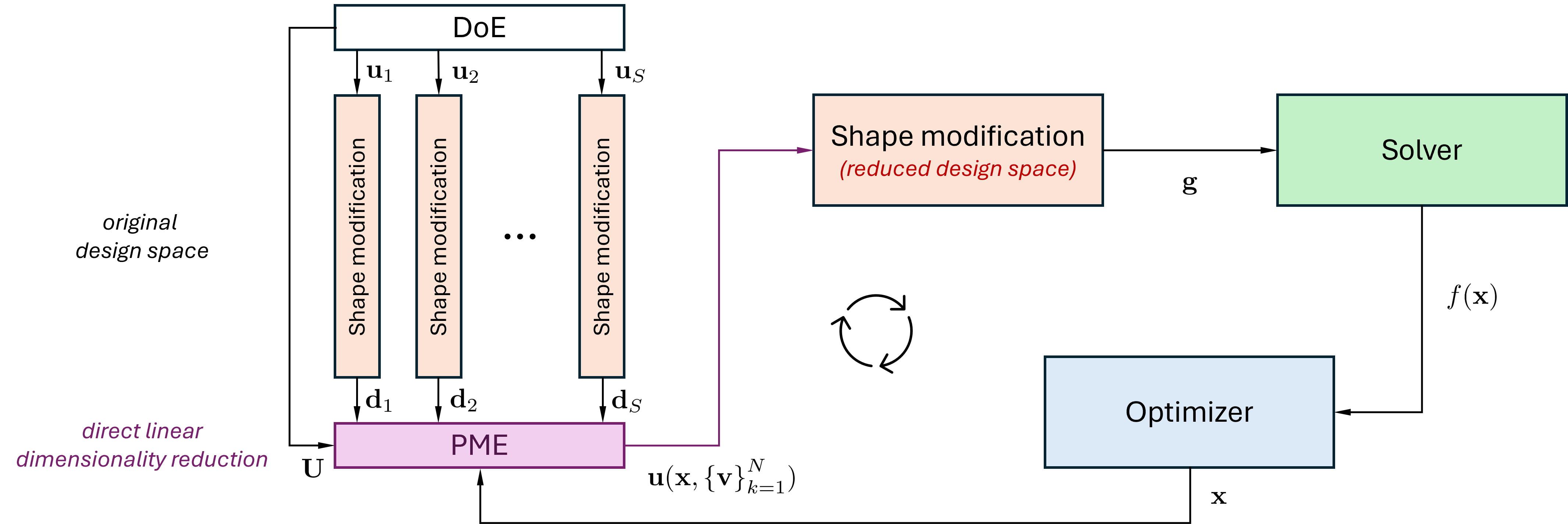}
    \caption{\change{Block diagram for the solution of shape optimization problems through PME}}
    \label{fig:block_pme}
\end{figure*}

Starting from the discrete formulation detailed in the preceding section, the dimensionality reduction procedure employs the data matrix $\mathbf{D}$ augmented with values of design parameters from the original parameterization. 



Practically, defining $\hat{\bfu}=\bfu-\langle\bfu\rangle$ (as already done for the shape modification vector), the embedding is achieved by introducing a new matrix $\mathbf{P}$ of dimensionality $\left[(3L+M)\times S\right]$ as follows
\begin{equation}
\mathbf{P}=\left[
\begin{array}{c}
\mathbf{D}\\
\mathbf{U}
\end{array}
\right]
\qquad
\mathrm{with}
\qquad
\mathbf{U}=\left[
\begin{array}{ccc}
\hat{u}_{1,1} & & \hat{u}_{1,S}\\
\vdots & \cdots & \vdots\\
\hat{u}_{M,1} & & \hat{u}_{M,S}\\
\end{array}
\right]
\end{equation}   
where the matrix $\mathbf{U}$ of the original design variables is added to the data matrix $\mathbf{D}$ with a null weight $\mathbf{W}_\bfu$ such that
\begin{equation}\label{eq:pme_weight}
\mathbf{W}_\bfu=\mathbf{0}
\qquad
\mathrm{and}
\qquad
\widetilde{\mathbf{W}}=\left[
\begin{array}{cc}
\mathbf{W} & \mathbf{0}\\
\mathbf{0} & \mathbf{W}_\bfu\\
\end{array}
\right]
\end{equation}   
and so recasting Eq. \ref{eq:pca} to
\begin{equation}\label{eq:pme_pca}
\widetilde{\mathbf{A}}\widetilde{\mathbf{G}}\widetilde{\mathbf{W}}\widetilde{\mathbf{Z}}=\widetilde{\mathbf{Z}}\widetilde{\boldsymbol{\Lambda}} \qquad \mathrm{with} \qquad \widetilde{\mathbf{A}}=\frac{1}{S}\mathbf{PP}^\mathsf{T} 
\end{equation} 
where
\begin{equation}\label{eq:pme}
\widetilde{\mathbf{G}}=\left[
\begin{array}{cc}
\mathbf{G} & \mathbf{0}\\
\mathbf{0} & \mathbf{I}\\
\end{array}
\right]
\qquad
\mathrm{and}
\qquad
\widetilde{\mathbf{Z}} = \left[\tilde{\bfz}_1 \,\,\, \dots \,\,\, \tilde{\bfz}_S \right]
\end{equation} 
with 
\begin{equation}
\tilde{\bfz}_k=\left[
\begin{array}{c}
\bfz_k\\
\bfv_k
\end{array}
\right]
\end{equation} 

Having given a null weight to $\mathbf{U}$ does not remove the contribution of the design variables from the inner product, but just cancels as many columns as $M$ from the matrix $\widetilde{\mathbf{A}}\widetilde{\mathbf{G}}\widetilde{\mathbf{W}}$, thus Eqs. \ref{eq:pca} and \ref{eq:pme} provides the same eigenvalues ($\boldsymbol{\Lambda}=\widetilde{\boldsymbol{\Lambda}}$) and geometrical components of the eigenvectors ($\bfz_k$, except for a multiplicative constant). The proof of the equivalence between KLE/PCA and PME is given in \cite{serani2023parametric}. In addition and as desired, the solution of Eq. \ref{eq:pme} provides the eigenvector components $\bfv_k$ that embed the original design variables $\bfu$. 

\begin{figure*}[!b]
    \centering
    \includegraphics[width=0.32\textwidth]{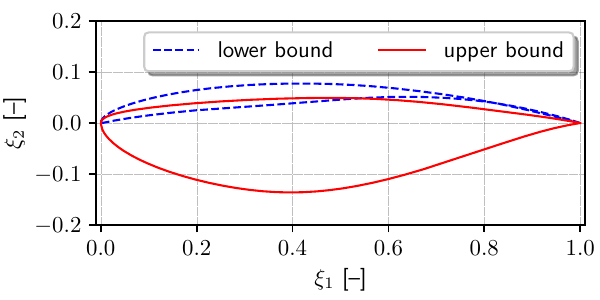}
    \includegraphics[width=0.32\textwidth]{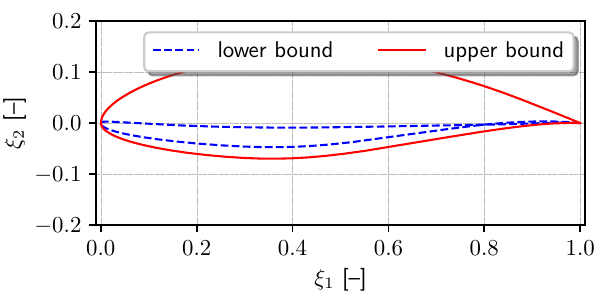}
    \includegraphics[width=0.32\textwidth]{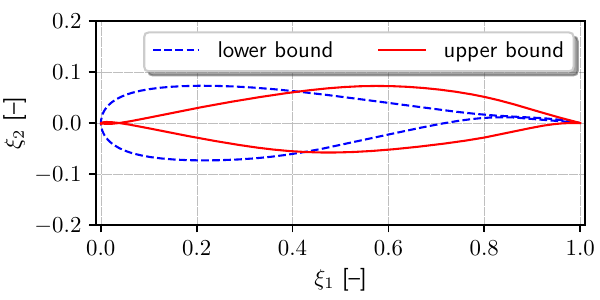}\\
    \includegraphics[width=0.32\textwidth]{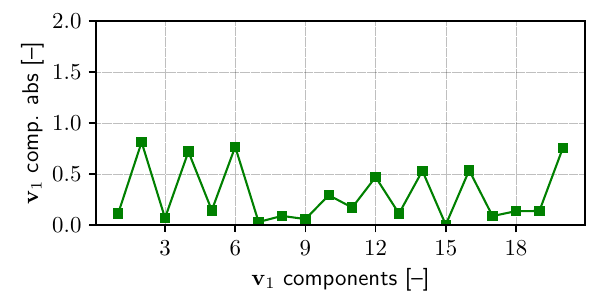}
    \includegraphics[width=0.32\textwidth]{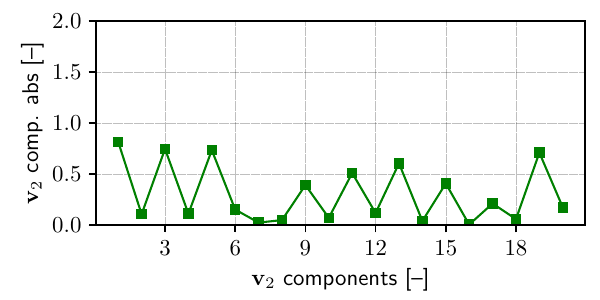}
    \includegraphics[width=0.32\textwidth]{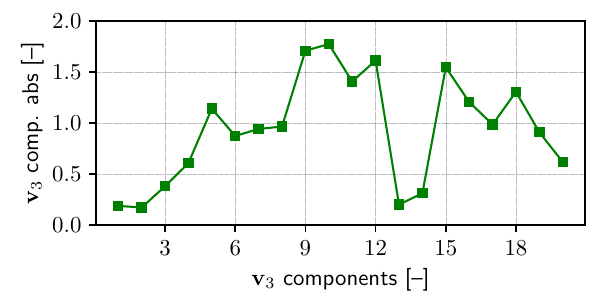}\\
    \caption{Comparison of (top) PCA geometric and (bottom) PME design variables eigenvector components \cite{serani2023efficient}}
    \label{fig:pme}
\end{figure*}
To reconstruct at least all the samples in $\mathbf{D}$, the coefficients $\boldsymbol{\theta}_j$, for $j=1, \dots, S$, are evaluated projecting the matrix $\mathbf{P}$ on $\widetilde{\mathbf{Z}}'$, that contains only the first $N$ eigenvectors of $\widetilde{\mathbf{Z}}$, retaining the desired level of variance of the original design space, as follows
\begin{equation}
\boldsymbol{\Theta}=\mathbf{P}^\mathsf{T}\widetilde{\mathbf{G}}\widetilde{\mathbf{W}}\widetilde{\mathbf{Z}}' 
\end{equation}
with $\boldsymbol{\Theta}=[\boldsymbol{\theta}_1 \, \dots \, \boldsymbol{\theta}_S]^\mathsf{T}$. Consequently, the reduced design variables $\mathbf{x} = [x_1 \, \dots \, x_N]^\mathsf{T}$ can be bounded such as 
\begin{equation}
\underset{j}{\min} \, \Theta_{jk} \leq x_k \leq \underset{j}{\max} \, \Theta_{jk} \qquad k=1, \, \dots \, N.
\end{equation}

The PME of the original design variables is finally achieved by
\begin{equation}\label{eq:recu}
\bfu \approx \check{\bfu}=\langle{\bfu}\rangle+\sum_{k=1}^N x_k\bfv_k
\end{equation}
where the components of the eigenvector $\bfv_k$ embed (or contain) the reduced-order representation of the original design parameterization. \change{Figure \ref{fig:block_pme} shows the integration of the original design variables in the dimensionality reduction procedure and how the optimization loop is modified by retrieving the original parametrization through PME eingevectors}.

The method has been presented in \cite{serani2023parametric} and successfully applied for both the aerodynamics optimization of NACA and RAE \cite{serani2024aerodynamic} foils and the hydrodynamic optimization of a destroyer-type vessel \cite{serani2024hydrodynamic}. \change{Figure \ref{fig:pme} provide an example that compares PCA and PME methods by visualizing the first three eigenvector components for both approaches applied to the shape optimization of the RAE-2822 airfoil \cite{serani2024aerodynamic}. In particular, the top row shows the shape modification associated with the PCA re-parameterization. The bottom row, instead, gives the eigenvectors $\bfv_k$ that enable the explicit mapping between reduced coordinates and original design variables. Furthermore, looking at peaks of each $\bfv_k$ components, PME allows to assess which of the original design parameters mostly influence the associated portion of geometric variance retained.}

\begin{figure*}[!b]
    \centering
    \includegraphics[width=1\textwidth]{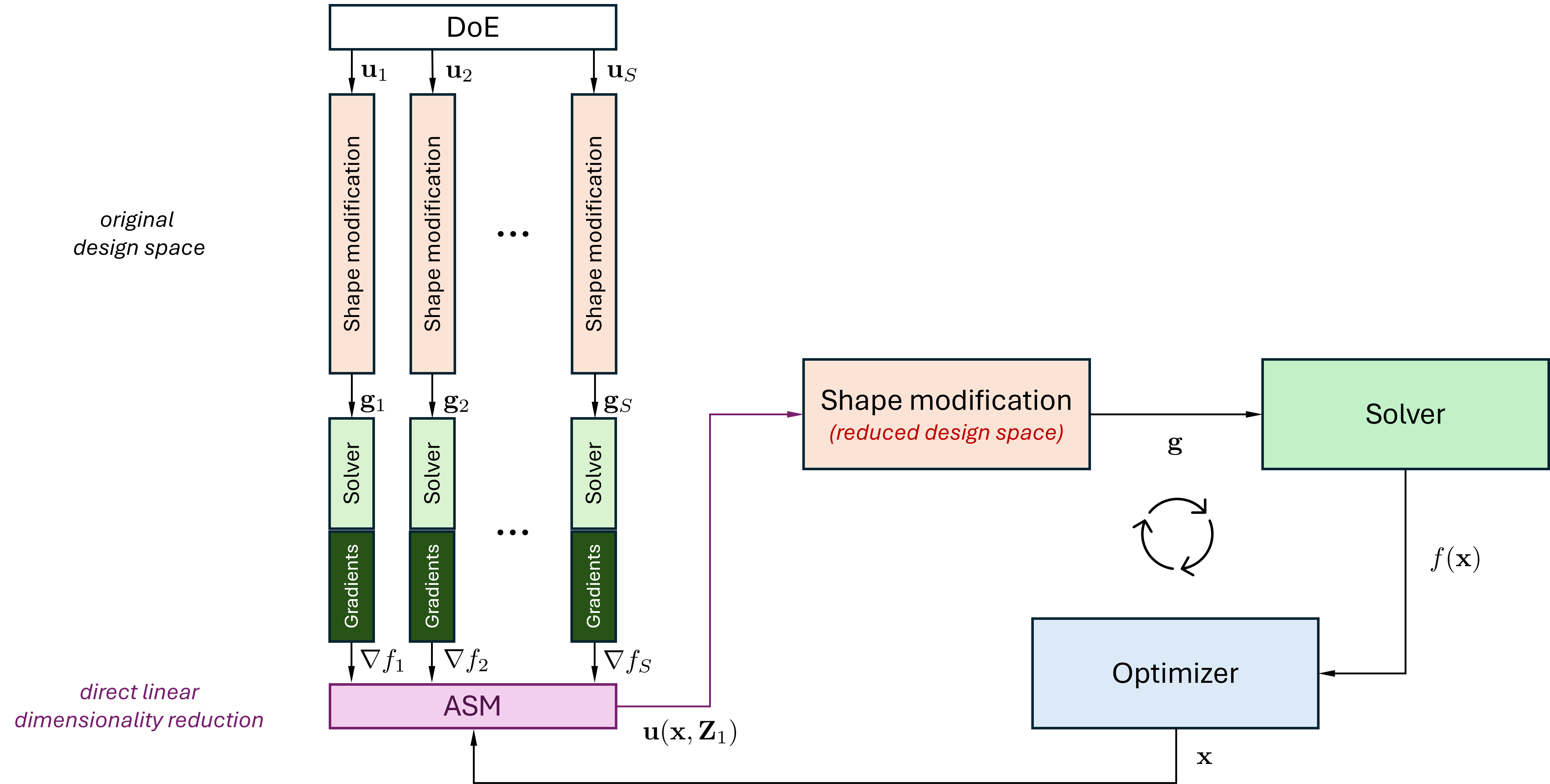}
    \caption{\change{Block diagram for the solution of shape optimization problems through ASM}}
    \label{fig:block_asm}
\end{figure*}
\subsubsection{Active subspace method}

The active subspace method (ASM) is a powerful dimensionality reduction technique tailored for design optimization \cite{constantine2014active}, especially in engineering applications where the relationship between input design variables and a scalar output of interest (e.g., performance metric) is governed by complex, possibly nonlinear functions. The primary goal of ASM is to identify directions in the input space along which the output function varies most significantly. These directions, known as active subspaces, enable the reduction of the problem's dimensionality by projecting the high-dimensional input space onto a lower-dimensional subspace without significant loss of critical information about the system's behavior.

Consider the objective function $f$ of problem \ref{eq:soform} that models the system of interest. Assuming $f$ as differentiable with respect to the inputs $\bfu$, the gradient of $f$, denoted by
\begin{equation}
\nabla f(\bfu)= \left[\frac{\partial f}{\partial u_1} \,\, \dots \,\, \frac{\partial f}{\partial u_M}\right]^\mathsf{T},    
\end{equation}
contains information about the sensitivity of the output with respect to changes in each input variable.

Similarly to PCA, the core idea of ASM revolves around the analysis of the covariance matrix of the gradient, defined as:
\begin{eqnarray}\label{eq:ASM}
\mathbf{C}=\mathbb{E}\left[\nabla f(\bfu)\nabla f(\bfu)^\mathsf{T}\right],   
\end{eqnarray}
where $\mathbb{E}[\cdot]$ denotes the expected value, assuming a probability distribution for the inputs $\bfu$ and consequently can be estimated through MC sampling as
\begin{equation}
\mathbf{C} \approx \widetilde{\mathbf{C}} = \frac{1}{P}\sum_{j=1}^P \left(\nabla f(\bfu_j)\nabla f(\bfu_j)^\mathsf{T}\right)   
\end{equation}
with $P$ the number of MC samples.

To identify the active subspace, the eigenvalue decomposition of the approximated covariance matrix $\widetilde{\mathbf{C}}$ is performed as follows:
\begin{equation}\label{eq:ASM}
\widetilde{\mathbf{C}}=\mathbf{Z}\mathbf{\Lambda}\mathbf{Z}^\mathsf{T}    
\end{equation}
where $\mathbf{Z}$ is the matrix whose columns are the eigenvectors of $\widetilde{\mathbf{C}}$, and $\mathbf{\Lambda}$ is a diagonal matrix containing the corresponding eigenvalues $\lambda_i$, sorted in descending order. The eigenvectors associated with the largest eigenvalues represent the directions in the input space that have the most significant impact on the variation of the output function $f$.

Partitioning the eigenspace into active and inactive subspaces, the active subspace is spanned by the eigenvectors corresponding to the largest eigenvalues, which capture the most significant directions of variation; in contrast, the inactive subspace is spanned by the remaining eigenvectors. Mathematically, the active subspace is defined by selecting the first $N<M$ eigenvectors, forming the matrix $\mathbf{Z}_1\in\mathbb{R}^{M\times N}$, and similarly, the inactive subspace is represented by $\mathbf{Z}_2\in\mathbb{R}^{M\times (M-N)}$.

\begin{figure*}[!b]
    \centering
    \includegraphics[width=0.75\textwidth]{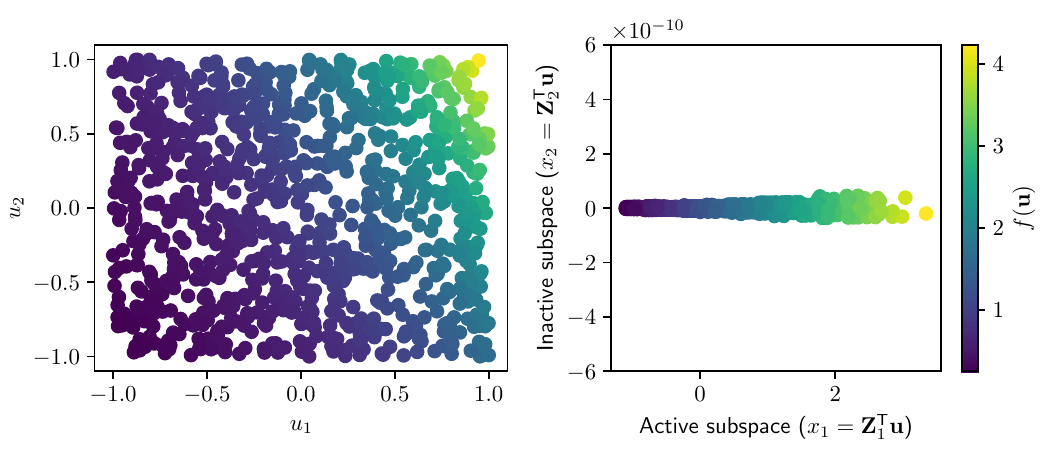}
    \caption{Example of ASM application on bi-dimensional problem}
    \label{fig:asm}
\end{figure*}
Once the active subspace is identified, the original high-dimensional optimization problem can be reformulated in terms of the reduced variables 
\begin{equation}\label{eq:backasm}
\bfx = \mathbf{Z}_1^\mathsf{T}\bfu    
\end{equation}
This transformation enables the exploration of the most influential directions with fewer variables, significantly reducing the computational complexity of the optimization process.  

\change{Figure \ref{fig:block_asm} shows the block diagram associated with ASM dimensionality reduction procedure and subsequent optimization loop, where the original parameterization is retrieved through back-mapping as per Eq. \ref{eq:backasm}. 
A simple example illustrating how the ASM operates on a bi-dimensional design space is shown in Fig. \ref{fig:asm}. The figure clearly shows that the gradient of the objective function predominantly follows the diagonal direction of the design space. This indicates that the function can be effectively represented in a one-dimensional active subspace, defined by a single reduced variable.}

The ASM offers a rigorous framework for identifying the most important directions in the design space, enabling efficient optimization by focusing computational efforts on the subspace that most significantly affects the output. This approach is particularly beneficial in shape optimization, where the relationship between design variables and performance metrics can be highly complex and nonlinear. By reducing the problem's dimensionality, ASM facilitates a more manageable and computationally efficient optimization process. 

It may be emphasized that the ASM relies on the computation of gradients of the objective function with respect to the design variables. To accurately estimate the gradient covariance matrix, a representative set of samples from the design space must be evaluated. The number of samples is crucial for capturing the variability and sensitivity of the objective function across the design space. However, the requirement for a sufficient number of samples introduces challenges, particularly in high-dimensional spaces where the number of required evaluations grows.
The estimation of gradients at each sample point is computationally intensive, particularly in high-dimensional settings where the curse of dimensionality comes into play. For a design space with $M$ dimensions, evaluating the gradient at a single point typically requires $M+1$ function evaluations (assuming finite differences are used for gradient approximation), which quickly becomes prohibitive as $M$ increases. This challenge necessitates efficient strategies for sampling and gradient estimation to make ASM practical for real-world optimization problems, such as adaptive sampling, surrogate modeling, and sparse gradient and regularization. Adaptive sampling techniques can be employed to strategically select sample points where the objective function's behavior is most uncertain or where gradients are likely to provide the most information. This approach aims to minimize the number of required function evaluations while still capturing the essential features of the design space. Surrogate models can approximate the objective function and its gradients. Once trained on a set of initial samples, these models can provide gradient estimates at additional points without direct evaluation of the expensive objective function, thereby reducing computational costs. In some applications, the objective function may have sparse gradients, meaning that only a few directions significantly influence the output. Identifying and focusing on these directions can reduce the effective dimensionality. Regularization techniques can also be employed to stabilize the estimation of the covariance matrix when the number of samples is limited.

The need for comprehensive sampling and gradient estimation underscores the trade-off between the accuracy of the active subspace identification and computational efficiency. While ASM offers a powerful framework for reducing dimensionality and facilitating optimization, its practical application must consider these sampling and gradient estimation challenges. By adopting strategies such as adaptive sampling, surrogate modeling, and focusing on sparse gradients, the computational burden of ASM can be managed, making it a viable option for optimizing complex engineering designs.
While ASM provides a methodological advancement in design-space dimensionality reduction, its application requires careful consideration of the preliminary design-space assessment. Addressing the challenges associated with sampling and gradient estimation is essential for harnessing the full potential of ASM in shape optimization.

ASM has been effectively used to reduce high-dimensional design spaces in external and internal \cite{seshadri2020supporting,lopez2022global} aerodynamic optimizations, facilitating the exploration and improvement of airfoil \cite{demo2021supervised,alswaitti2022dimensionality} and wing \cite{grey2018active, li2019surrogate, mufti2024multifidelity}, aero-engine nacelle \cite{berguin2015dimensionality, tejero2022aerodynamic}, as well as satellite applications \cite{hu2018conceptual}, including also fluid-structure interaction \cite{boncoraglio2021model} and uncertainty quantification \cite{song2024sensitivity} problems. This has helped in achieving better performance metrics like lift-to-drag ratios by focusing on the most influential variables, leading to computationally efficient designs.
In naval engineering, ASM has supported the optimization of ship hulls and other marine structures, simplifying the complex sets of design variables into manageable subspaces and thus improving hydrodynamic performance and structural integrity \cite{demo2021hull,tezzele2023multifidelity}. \change{Additionally, in the field of internal combustion engines, ASM has recently been employed to analyze and reduce the dimensionality of complex design problems related to the combustion chamber structure of Wankel rotary engines, facilitating the identification of critical structural parameters and significantly enhancing optimization efficiency \cite{li2024research}.}

These examples collectively underscore ASM's capability in reducing computational expenses, simplifying complex optimization problems, and yielding improved design efficiencies across diverse engineering disciplines.


\subsection{Nonlinear dimensionality reduction methods}
Nonlinear dimensionality reduction techniques offer enhanced flexibility for capturing complex relationships in shape data. Examples include local PCA and kernel PCA, which extend linear methods to capture nonlinear shape variations more effectively, as well as neural network generalization of PCA through autoencoders, a class of neural networks that learn hierarchical representations of shapes, enabling nonlinear dimensionality reduction through successive encoding and decoding layers. The following subsection provides details of nonlinear dimensionality reduction methods used in shape optimization literature.

\begin{figure*}[!t]
    \centering
    \begin{minipage}[!t]{0.48\textwidth}
    \includegraphics[width=1\textwidth]{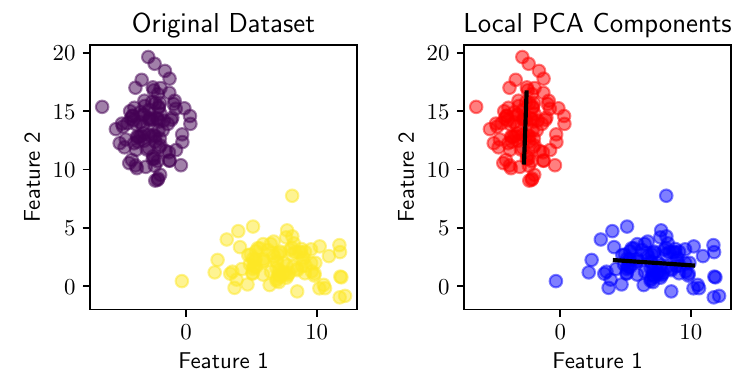}
    \caption{LPCA examples with first components}
    \label{fig:lpca}
    \end{minipage}
    \hfill
    \begin{minipage}[!t]{0.48\textwidth}
    \includegraphics[width=1\textwidth]{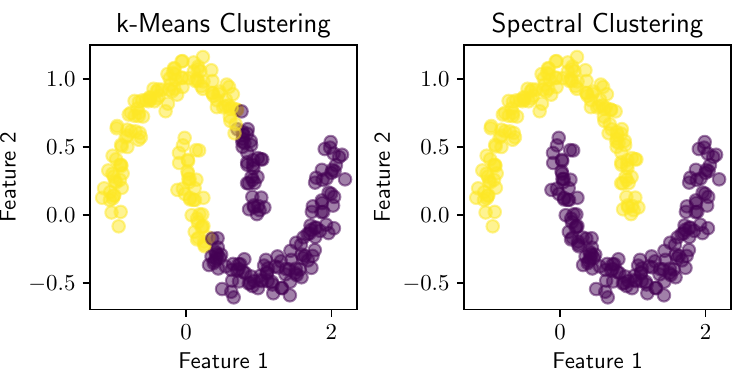}
    \caption{Different clustering approaches}
    \label{fig:cluster}
    \end{minipage}
    \hfill
\end{figure*}
\begin{figure*}[!b]
    \centering
    \includegraphics[width=1\textwidth]{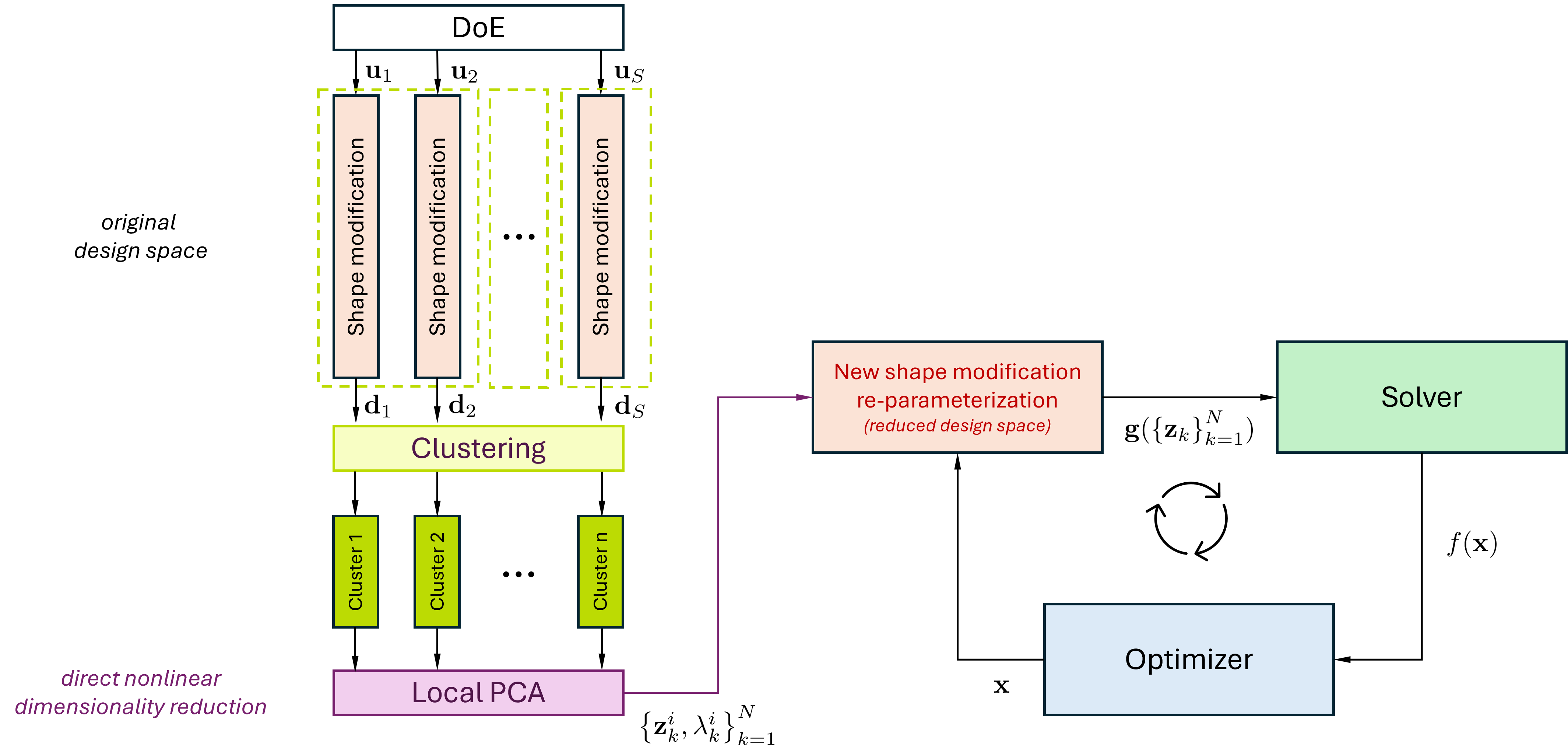}
    \caption{\change{Block diagram for the solution of shape optimization problems through LPCA}}
    \label{fig:block_lpca}
\end{figure*}
\subsubsection{Local principal component analysis}
Local PCA (LPCA) performs a PCA for every different disjoint region of the input space, assuming that, if the local regions are small enough, the data manifold will not curve much over the extent of the region and the linear model will be a good fit \cite{kambhatla1997dimension}. 

The first step in LPCA is to cluster the data in $k$ sets, applying a clustering algorithm, such that $\mathbf{D} = \{\mathbf{D}_1,\dots,\mathbf{D}_i\}_{i=1}^k$, as shown in Fig. \ref{fig:lpca}. 

LPCA can be performed with different clustering techniques, such as the $k$-means algorithm \cite{lloyd1982least} or a spectral clustering \cite{ng2001spectral}. 
One issue in $k$-means is that using the Euclidean distance as a similarity measure assumes a convex shape to the underlying clusters \cite{aggarwal2013data}, which cannot be guaranteed a priori. On the contrary, spectral clustering can be effective even if the clusters' shapes are more complex, see Fig. \ref{fig:cluster}. There are several versions of the spectral clustering algorithms, the main difference is in which graph Laplacian is used \cite{von2007tutorial}, e.g., the symmetric normalized Laplacian \cite{ng2001spectral}.

After the data are partitioned in $k$ clusters, a PCA is performed on them, as per Eq. \ref{eq:pca}, solving $k$ PCA eigenproblems
\begin{equation}
\mathbf{A}_i\mathbf{G}_i\mathbf{W}_i\mathbf{Z}_i= \mathbf{Z}_i\boldsymbol{\Lambda}_i  \qquad \forall i= 1,\dots ,k
\end{equation}

LPCA is particularly useful in shape optimization where the design space is inherently complex and multidimensional. By localizing the PCA, LPCA can adapt more flexibly to local variations within the data, capturing essential features that global PCA might miss. Each local model can offer insights into specific regions of the design space, potentially revealing localized trends and dependencies that are relevant for detailed design analysis. \change{Figure \ref{fig:block_lpca} depicts the block diagram using LPCA, involving partitioning the design space into distinct clusters and performing PCA locally within each cluster.}

While LPCA offers several benefits, it comes with challenges that need careful consideration. The effectiveness of LPCA heavily depends on the initial clustering. Poor clustering can lead to misleading results and inefficient dimensionality reduction. Determining the optimal number of clusters $k$ is crucial. Too many clusters can lead to overfitting and increased computational complexity, while too few may not adequately capture the necessary local variations. Although LPCA can reduce the complexity within each local model, managing multiple local analyses simultaneously increases the overall computational burden, especially when the number of clusters is large.

\begin{figure*}[!b]
    \centering
    \includegraphics[width=1\textwidth]{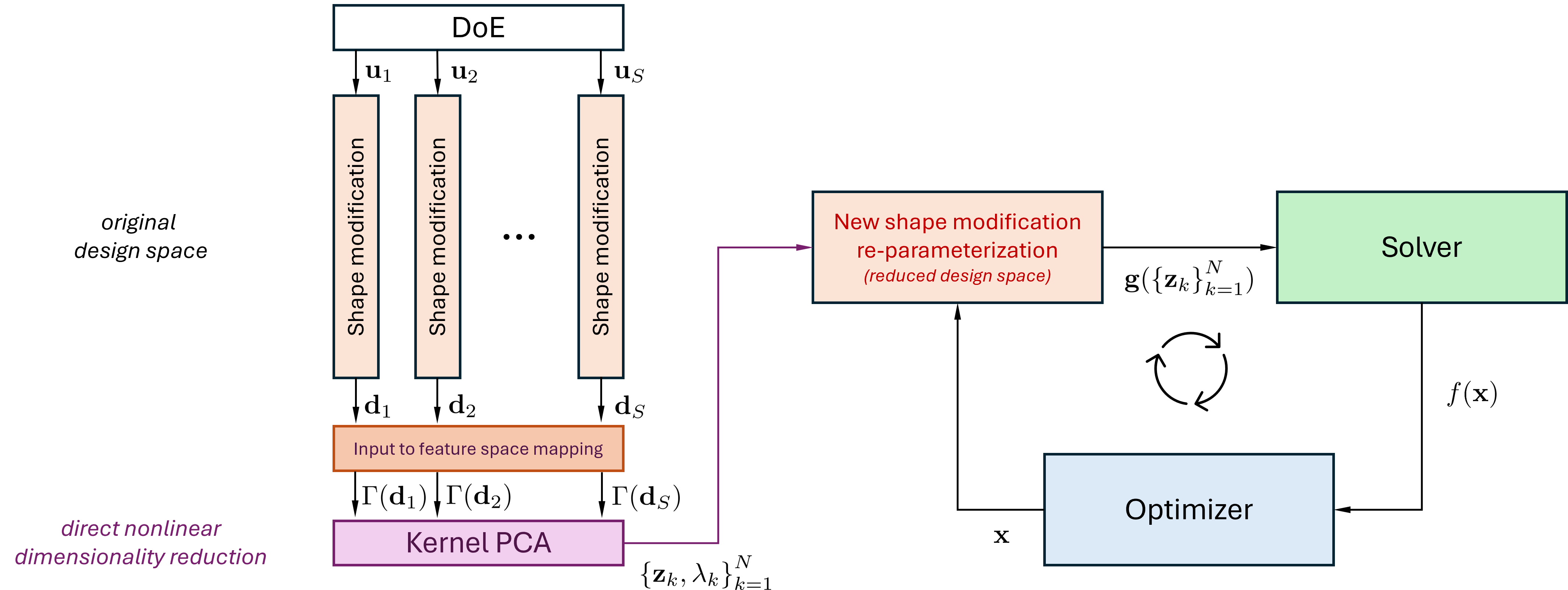}
    \caption{\change{Block diagram for the solution of shape optimization problems through KPCA}}
    \label{fig:block_kpca}
\end{figure*}
LPCA has been successfully applied as a preliminary tool for design-space assessment and dimensionality reduction in \cite{d2018nonlinear} for marine engineering applications. Nevertheless, the results are highly dependent on the clustering procedure and especially on the number of clusters used. Moreover, the number of clusters $k$ should be set carefully to avoid extensive computation.

\subsubsection{Kernel principal component analysis}
The kernel PCA (KPCA) method \cite{scholkopf1998nonlinear} finds directions of maximum variance in a higher (possibly infinite) dimensional feature space $\mathcal{F}$, mapping the points from the input space $\mathcal{I}$ by a possible nonlinear function $\Gamma: \mathcal{I} \rightarrow \mathcal{F}$ as
\begin{equation}\label{e:map_kernel}
\mathbf{d}(\bfu_k) \rightarrow \Gamma(\mathbf{d}(\bfu_k)), ~~~~~~ \forall \, k = 1,\dots,S  
\end{equation}
Then PCA is computed in the feature space $\mathcal{F}$ \change{(see Fig. \ref{fig:block_kpca})}.
Assuming that $\sum_k \Gamma(\mathbf{d}(\bfu_k))=0$, the kernel principal component $[\mathbf{z}_1 \, \dots \, \bfz_J]$ can be find solving the eigenproblem
\begin{equation}\label{e:eigen1_KPCA}
\mathbf{\Sigma}_\Gamma \mathbf{z}_j = \lambda_j \mathbf{z}_j 
\end{equation}
where $\mathbf{\Sigma}_\Gamma$ is the $[J\times J]$ covariance matrix in the feature space $\mathcal{F}$, defined as
\begin{equation}\label{e:cov_kernel}
\mathbf{\Sigma}_\Gamma = \frac{1}{S}\sum_{k=1}^{S}\Gamma(\mathbf{d}(\bfu_k))\Gamma(\mathbf{d}(\bfu_k))^\mathsf{T}
\end{equation}

KPCA allows the solution of Eq. \ref{e:eigen1_KPCA} without computing explicitly the Eq. \ref{e:map_kernel}, since it appears only within an inner product \cite{smola1998learning}, which can be computed efficiently by a kernel function
%
$K(\mathbf{d}(\bfu_i),\mathbf{d}(\bfu_k)) = \Gamma(\mathbf{d}(\bfu_i))^\mathsf{T}\Gamma(\mathbf{d}(\bfu_k))$.
%
Defining $\mathbf{z}_j$ as a linear expansion of $\Gamma(\mathbf{d}(\bfu_k))$ 
\begin{equation}\label{e:lin_comb_u_m}
\mathbf{z}_j = \sum_{k=1}^{S}c_{jk}\Gamma(\mathbf{d}(\bfu_k))
\end{equation}
the Eq. \ref{e:eigen1_KPCA} can be recasted as
\begin{equation}\label{e:eigen2_KPCA}
\mathbf{K}\mathbf{c}_j = \lambda_j S \mathbf{c}_j
\end{equation}
where $\mathbf{K}$ is the symmetric and positive-semidefinite $[S\times S]$ kernel matrix, with  $\mathbf{K}_{ik} =  K(\mathbf{d}(\bfu_i),\mathbf{d}(\bfu_k))$. The length of the $S$-component vector $\mathbf{c}_j$ is chosen such that $\mathbf{z}_j^\mathsf{T}\mathbf{z}_j=\lambda_j S \mathbf{c}_j^\mathsf{T}\mathbf{c}_j=1$. Once the eigenproblem in Eq. \ref{e:eigen2_KPCA} is solved, the new design variables can be found projecting $\Gamma(\mathbf{d}(\bfu))$ on $\mathbf{z}_j$ as
\begin{align}\label{e:key}
\bfx & = \Gamma(\mathbf{d}(\bfu))\mathbf{z}_j = \sum_{k=1}^{S}c_{jk}\Gamma(\mathbf{d}(\bfu)^\mathsf{T}\Gamma(\mathbf{d}(\bfu_k)) \nonumber \\
& =\sum_{k=1}^{S}c_{jk} K(\mathbf{d}(\bfu),\mathbf{d}(\bfu_k)) 
\end{align}
\begin{figure*}[!b]
    \centering
    \includegraphics[width=0.75\textwidth]{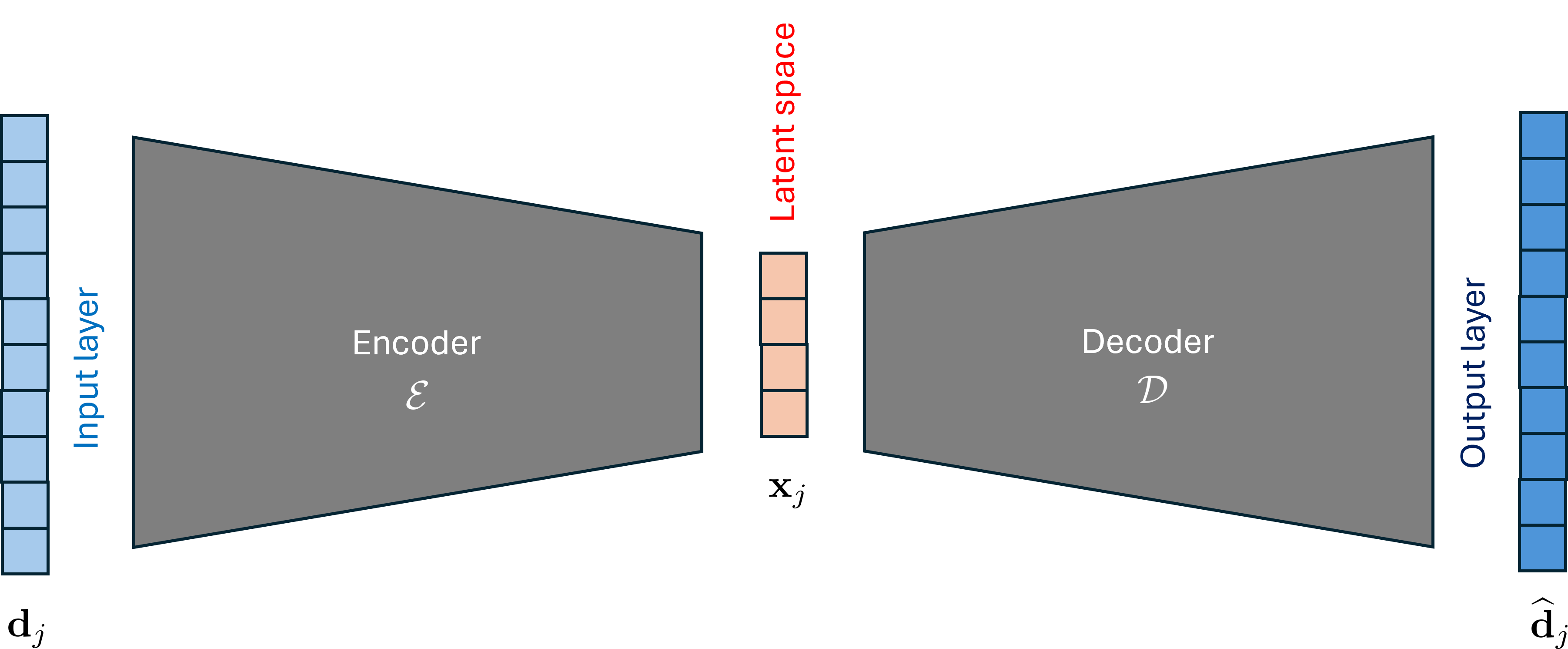}
    \caption{Example scheme of autoencoder architecture}
    \label{fig:dae}
\end{figure*}
The reconstruction of the original data from the feature space $\mathcal{F}$ in KPCA is more problematic than PCA, since it needs to find, for every point $\Gamma(\mathbf{d}(\bfu_k))$, the relative pre-image $\mathbf{d}(\bfu_k)$ in the input space $\mathcal{I}$.  
Nevertheless, KPCA has been successfully applied for the design-space assessment of a destroyer-type vessel in \cite{d2018nonlinear}, where the approximate pre-images technique proposed in \cite{bakir2004learning} is used, and aerodynamic optimization of NACA foil \cite{gaudrie2020modeling}, \change{RAE-2822 and CRM wing \cite{zhao2024supervised}}.

\begin{figure*}[!t]
    \centering
    \includegraphics[width=1\textwidth]{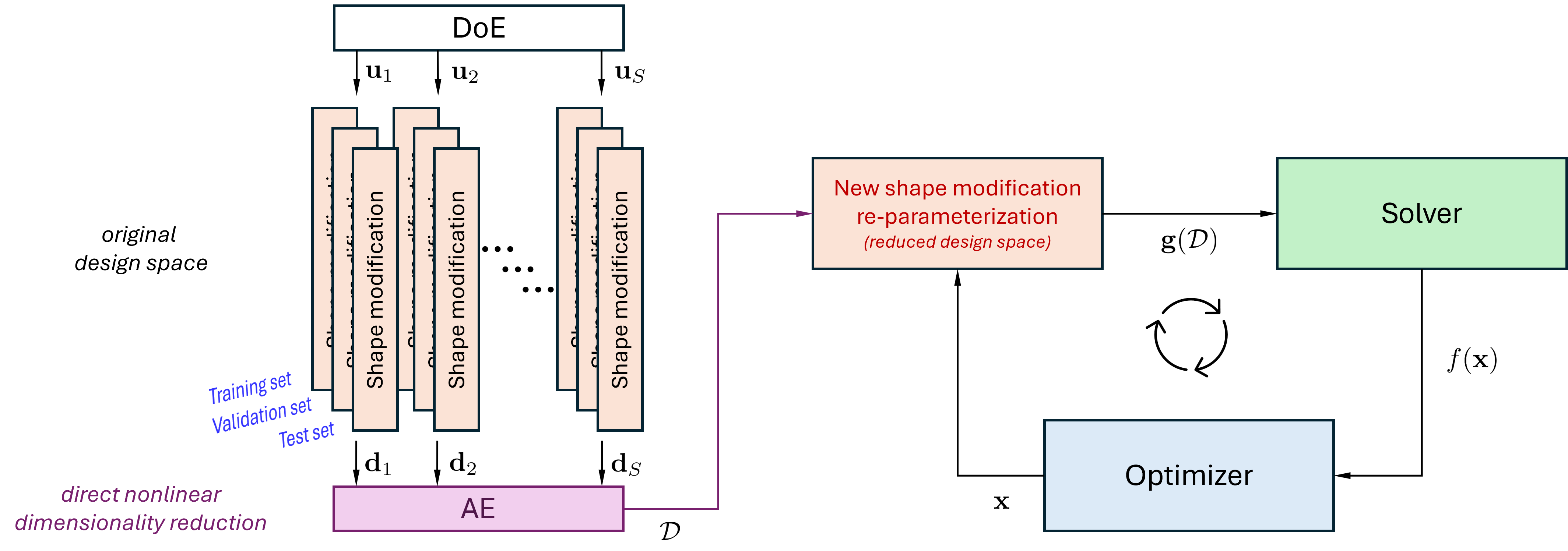}
    \caption{\change{Block diagram for the solution of shape optimization problems through AE}}
    \label{fig:block_ae}
\end{figure*}
\subsubsection{Autoencoders}
Autoencoders (AEs) are a type of artificial neural networks that are used primarily for the task of data compression and decompression through a dual-component architecture \cite{HintonSalakhutdinov2006b} (see Fig. \ref{fig:dae}). This architecture consists of two main functions: (1) an encoder function $\mathcal{E}$ maps the data $\mathbf{d}$ to compress data $\bfx$; (2) a decoder function $\mathcal{D}$ maps from the compressed data $\bfx$ back to $\hat{\mathbf{d}}$. This operation is performed by setting the same number of neurons in the input and output layer and constraining the hidden layer to have $N<M$ neurons.

Consider a single hidden layer AE, if the new design variable $\bfx$ can be written as 
\begin{equation}
\bfx = \mathcal{E}(\mathbf{H}^{(1)}\mathbf{d} + \mathbf{b}^{(1)})
\end{equation}
where $\mathbf{H}$ is a relative weight matrix, $\mathbf{b}$ the bias vector, and the apex ``(1)'' represent the hidden layer, then the reconstruction vector $\hat{\mathbf{d}}$ from $\bfx$ can be expressed as 
\begin{equation}
\hat{\mathbf{d}} = \mathcal{D}(\mathbf{H}^{(2)}\bfx + \mathbf{b}^{(2)})
\end{equation}
where the apex ``(2)'' represent the output layer. 
The network parameters $\mathbf{H}$ and $\mathbf{b}$, are evaluated minimizing the reconstruction error

\begin{align}\label{e:rec_err_DAE}
E(\mathbf{H}^{(1)},\mathbf{b}^{(1)},\mathbf{H}^{(2)},\mathbf{b}^{(2)}) 
= \frac{1}{2}\sum_{k=1}^{S}\|\mathbf{d}_k - \hat{\mathbf{d}}_k\|^2  
\nonumber\\
= \frac{1}{2}\sum_{k=1}^{S}\|\mathbf{d}_k - \mathcal{D}(\mathbf{H}^{(2)}\mathcal{E}(\mathbf{H}^{(1)}\mathbf{d}_k + \mathbf{b}^{(1)}) + \mathbf{b}^{(2)})\|^2 
\end{align}

If $\mathcal{E}$ and $\mathcal{D}$ are linear then the Eq. \ref{e:rec_err_DAE} has a unique global minimum, in which the weights in the hidden layer span the same subspace as the first $N$-principal components of the data \cite{bourlard1988auto,baldi1989neural}. AE with nonlinear activation functions and more hidden layers (called deep autoencoder, DAE) provides a nonlinear generalization of the PCA \cite{lecun2015deep}, but in this case the error function (Eq. \ref{e:rec_err_DAE}) becomes nonconvex and the optimization algorithm may get stuck in poor local minima. 





\change{DAEs have been effectively employed for marine engineering problems involving single-hull \cite{d2018nonlinear,seo2024study} and multi-hull designs \cite{abbas2023deep}. In aerodynamic optimization contexts, DAEs have been successfully applied to the design of RAE 2822 airfoils \cite{yamazaki2023comparative}, transonic wing optimization \cite{karafi2024simultaneous}, aeroacoustic shape optimization of wind turbine blades \cite{kou2023aeroacoustic}, and rotorcraft airfoil aeroacoustic optimization tasks \cite{liu2024aerodynamic}. Recent studies have further advanced the field by integrating convolutional autoencoders for improved ship hull shape transformations \cite{seo2024study} and combining autoencoder architectures with physics-informed neural networks for efficient aerodynamic airfoil optimization \cite{liu2025cnn}, demonstrating their versatility and robustness across diverse optimization scenarios. Additionally, graph variational autoencoders have been recently employed for direct 3D aerodynamic shape optimization, leveraging their ability to efficiently represent complex, non-parametric 3D geometries and providing compact and smooth latent spaces suitable for advanced exploration via evolutionary algorithms \cite{jabon2024aerodynamic}.}

\begin{figure*}[!t]
    \centering
    \includegraphics[width=1\textwidth]{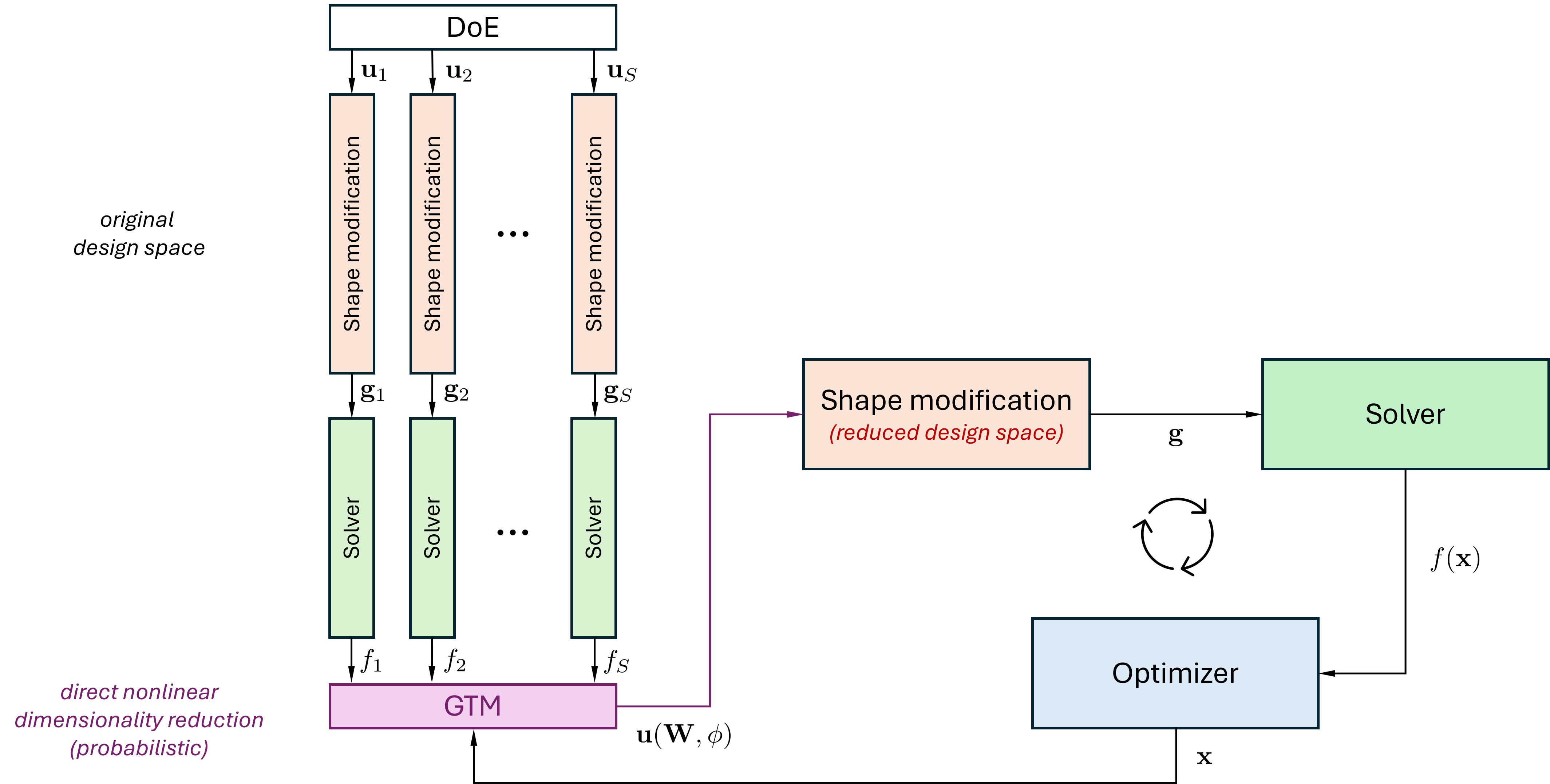}
    \caption{\change{Block diagram for the solution of shape optimization problems through GTM}}
    \label{fig:block_gtm}
\end{figure*}
Despite their versatility, AEs require careful tuning of parameters, including the number of neurons $N$ in the hidden layer, which corresponds to the intrinsic dimensionality of the data. The choice of this dimensionality, which is not known a priori, is crucial for balancing detail in the reconstructed output against the model's simplicity. Furthermore, the effectiveness of an AE can be significantly influenced by the initial weights, the configuration of the network layers\change{, and the training phase (see Fig. \ref{fig:block_ae})}. Continued advancements in AE technology, especially with the integration of regularization techniques and advanced training algorithms, promise to enhance their applicability in more complex and computationally demanding shape optimization tasks.

\subsubsection{Generative topographic mapping}
GTM is a sophisticated nonlinear dimensionality reduction technique that builds on the principles of probabilistic modeling to map high-dimensional data into a more manageable, lower-dimensional latent space. This approach is particularly relevant in the field of shape optimization, where each shape or design can be described by a complex set of parameters that define its geometry. GTM facilitates the exploration and optimization of these parameters by providing a systematic way to visualize and manipulate the data within a reduced space framework.

In the context of shape optimization, GTM operates by first defining a latent space, typically two or three-dimensional, which serves as the basis for data exploration. The high-dimensional data, representing different shape configurations, is assumed to be generated by a nonlinear transformation of this latent space. The transformation is mediated by a set of basis functions, often Gaussian, which map the grid points of the latent space to the original data space. The nonlinear mapping is expressed as 
\begin{equation}
\mathbf{y}(\mathbf{z}) = \mathbf{W} \phi(\mathbf{z}),    
\end{equation}
where $\mathbf{z}$ is a point in the latent space, $\mathbf{W}$ is a matrix of parameters, and $\phi(\mathbf{z})$ represents the basis functions.

This model not only helps in reducing the dimensionality but also in understanding the underlying structure of the shape space. By optimizing the parameters $\mathbf{W}$, GTM learns a manifold that best represents the distribution of the actual data points, thereby facilitating a smoother interpolation and visualization of the shape parameters’ effects. The fitting process involves maximizing the likelihood of the data under the model, typically by employing an expectation-maximization algorithm, which adjusts $\mathbf{W}$ and other model parameters to better fit the observed data \change{(see Fig. \ref{fig:block_gtm})}.

One of the key advantages of using GTM in shape optimization is its ability to provide insightful visualizations that help designers understand how different parameters influence the overall design. This is particularly useful in identifying regions of the parameter space that lead to optimal designs and in understanding how small variations in design parameters can lead to significant changes in the shape. Moreover, the probabilistic nature of GTM offers a way to estimate the density and distribution of different designs, providing a useful tool for stochastic optimization methods.

However, the application of GTM is not without challenges. The computational cost of setting up and training the GTM model can be substantial, especially as the number of dimensions and the complexity of the data increase. Additionally, the performance of GTM heavily relies on the correct setting of its hyperparameters, such as the number of grid points in the latent space and the variance of the basis functions. Inappropriate settings can lead to models that either overfit or fail to capture the essential variability of the data. Another potential drawback is that the nonlinear transformations inherent in GTM might introduce distortions that are difficult to interpret, particularly if the transformation hides important linear relationships or if the manifold learned by the model does not accurately reflect the true relationships in the higher-dimensional space. 

It has been applied for external \cite{viswanath2011dimension,zhao2024topographic} and internal \cite{viswanath2014constrained} aerodynamics problems, nevertheless, its success in practical applications hinges on the balance between model complexity and computational feasibility, making it essential for users to thoroughly evaluate the model’s behavior in the context of specific optimization tasks. \change{Recently \cite{zhang2024optimized}, an optimized GTM has been proposed to overcome some of these issues, introducing negative log-likelihood of data to measure the fitting quality between the nonlinear manifold with different dimensions and the high-dimensional space, as well as a variable-fidelity sample filtration to generate an efficient training set.}

In conclusion, while GTM presents a powerful tool for shape optimization by enabling effective dimensionality reduction and intuitive visualization of complex relationships, it requires careful management and understanding of its limitations.

\section{Discussion and Conclusions}
In the context of design-space dimensionality reduction for shape optimization problems, a general classification distinguishing between space and dimensionality reduction methods has been proposed. Traditional dimensionality reduction techniques, commonly utilized in statistics and machine learning, aim to reduce the number of variables (cardinality) while preserving as much information as possible from the original dataset. These methods, such as PCA, directly transform the design variable space (or input space) into a latent (lower-dimensional) space, through re-parameterization.
Conversely, space reduction techniques do not reduce the number of design variables. Instead, they constrain the range or domain within which these variables can vary. This approach simplifies the optimization problem by reducing the volume of the search space, thus reducing the problem's complexity by focusing the search and analysis efforts on a more manageable subset of the design space and potentially enhancing the efficiency of optimization algorithms without directly reducing the number of variables involved. For this reason, 
space reduction can be seen as an adjunct method that complements traditional dimensionality reduction techniques by narrowing the feasible design space, thereby indirectly contributing to the overall simplification and efficiency of the optimization process.

This survey underscores the role of dimensionality reduction in simplifying complex design spaces, thereby enabling efficient and effective shape optimization in designing functional surfaces. Traditional and emerging techniques, from indirect methods based on sensitivity analysis, passing through the most used direct approaches based linear methods like PCA and ASM, to advanced approaches like autoencoders, offer diverse tools for tackling the challenges associated with high-dimensional data. The integration of these techniques within simulation-based design optimization frameworks allows to alleviate the curse of dimensionality, improving optimization convergence, single- and multi-fidelity surrogate model training, as well as uncertainty quantification, enhancing both the speed and quality of the design process.

Linear dimensionality reduction methods, based on POD/KLE at the continuum, that are equivalent to PCA at the discrete level and can be solved through SVD, provide a solid foundation for identifying dominant modes of variation within a design space. However, they may not capture the nonlinearities inherent in complex shape optimizations. For this reason, nonlinear methods, such as local and kernel PCA, autoencoders, and generative topographic mapping, extend the capabilities of linear techniques, offering a more nuanced understanding of design spaces.

The novel parametric model embedding technique addresses a crucial gap in traditional PCA-based approaches by enabling a direct mapping from the reduced dimensionality space back to the original design variables. This advancement facilitates the use of parametric models within optimization workflows, enhancing the practical applicability of dimensionality reduction methods in industrial contexts. This pre-image problem can be solved also using the active subspace method, nevertheless, its computational complexity and cost are affected by the necessity of function values and associated gradients. 

Geometry-based PCA and its variants are found to be the most used, even for the implementation simplicity, nevertheless they can not be effective since no physical insight is taken into consideration in the reduced design space representation. To overcome this issue, further advancement has been achieved by integrating physical information within the design-space preliminary assessment. Physics-informed methods represent a promising frontier, incorporating domain-specific knowledge to guide the dimensionality reduction process. By integrating physical principles, these methods aim to ensure that the reduced design spaces remain meaningful and relevant to the optimization objectives \cite{serani2025_411}.

In conclusion, this survey illustrates the dynamic and multifaceted landscape of design-space dimensionality reduction techniques in shape optimization. Future research directions include the further development of physics-informed approaches, the exploration of novel nonlinear dimensionality reduction algorithms \cite{d2019augmented}, as well as the use of more sophisticated space reduction \cite{qiang2023multi}, and the integration of these methods into comprehensive optimization frameworks. \change{Promising extensions of traditional PCA to nonlinear settings, such as principal geodesic analysis, have recently emerged, employing advanced manifold learning methods like Grassmannian manifolds and separable shape tensors to capture complex geometric variations that linear PCA cannot adequately represent \cite{bo2024data,doronina2025aerodynamic}.} The continued evolution of these techniques will undoubtedly enhance their applicability and effectiveness in shape optimization and beyond, offering new opportunities for innovation in multidisciplinary optimization challenges, including also structural \cite{raghavan2012implicit,raghavan2013towards} and topology \cite{li2019dimension,munoz2022allying} optimization.

\section*{Acknowledgments}
The work was conducted in collaboration with the NATO-AVT-404 Research Task Group on ``Enhanced Design Processes of Military Vehicles through Machine Learning Methods''. 
The authors are grateful to the US Office of Naval Research
Global for its support through grants N62909-11-1-7011 and N62909-21-1-2042,
and acknowledge the support of the Italian Ministry of University and Research (MUR) through the PRIN 2022 BIODRONES  project, code 20227JNM52 (CUP: 853023005560006).

\bibliographystyle{unsrt}  
\bibliography{biblio}

\end{document}